\documentclass{artjlt}

\def\firstheadline{} 

\usepackage{amsfonts,amssymb,amsthm,amsmath}
\usepackage[matrix,arrow]{xy}

\theoremstyle{remark}

\def\ger{\mathfrak}
\DeclareMathOperator\pr{\mathrm{pr}}
\DeclareMathOperator\supp{\mathrm{supp}}

\makeatletter
\def\Set@Scallop[#1]#2#3{{#1}\Parens{#2}{#3}}
\newcommand\DeclareScalableOperator[2]{%
  \expandafter\def\csname#1\endcsname{\@ifnextchar[{{#2}\Set@Scallop}{{#2}\Set@Scallop[{}]}}
}
\makeatother

\DeclareScalableOperator{Ct}{\mathcal{C}} 
\DeclareScalableOperator{Cc}{\mathcal{C}_c} 
\DeclareScalableOperator{Hom}{\mathrm{Hom}} 
\DeclareScalableOperator{End}{\mathrm{End}} 
\DeclareScalableOperator{Aut}{\mathrm{Aut}} 
\DeclareScalableOperator{Der}{\mathrm{Der}} 
\DeclareScalableOperator{Ber}{\mathrm{Ber}} 

\newcommand{\Fa}{For all }
\newcommand{\fa}{for all }

\newcommand{\fs}{for some }
\newcommand\mathfa[1][{}]{\quad\text{\fa{#1} }}

\newcommand{\scth}{such that }
\newcommand{\AND}{and}

\newcommand\mathtxt[1]{\quad\text{{#1}}\quad}

\newcommand{\nd}{\mathtxt\AND}

\newcommand\vphi{\varphi}
\newcommand\vrho{\varrho}

\newcommand\eps{\varepsilon}

\newcommand\reals{\mathbb{R}}

\newcommand\vvoid{\varnothing}
\newcommand\sle{\leqslant}

\DeclareMathOperator\rk{\mathrm{rk}}
\DeclareMathOperator\Ad{\mathrm{Ad}}

\DeclareMathOperator\GL{\mathrm{GL}}
\DeclareMathOperator\id{\mathrm{id}}

\makeatletter
\newcommand\Size[7][1]{
                                 \ifx#20%
                                        \def\r@l{}\def\r@m{}\def\r@r{}%
                                 \else%
                                    \ifx#21%
                                           \def\r@l{\bigl}\def\r@r{\bigr}\def\r@m{\bigm}%
                                    \else%
                                           \ifx#22%
                                                 \def\r@l{\Bigl}\def\r@r{\Bigr}\def\r@m{\Bigm}%
                                            \else%
                                                 \ifx#23%
                                                        \def\r@l{\biggl}\def\r@r{\biggr}\def\r@m{\biggm}%
                                                  \else
                                                        \ifx#24%
                                                        \def\r@l{\Biggl}\def\r@r{\Biggr}\def\r@m{\Biggm}%
                                                        \fi%
                                                  \fi%
                                            \fi%
                                      \fi%
                                 \fi%
                                 \ifx#10%
                                       \def\r@m{}%
                                 \fi%
                                 \r@l#3{#4}\r@m#5{#6}\r@r#7%
}%

\makeatother

\newcommand\Set[3]{
                                 \Size{#1}{\{}{#2}{|}{#3}{\}}%
}%
\newcommand\Dual[3]{
                                 \Size[0]{#1}{\langle}{#2}{,}{#3}{\rangle}%
}%
\newcommand\Parens[2]{
  \Size[0]{#1}{(}{#2}{}{}{)}
}

\newcommand\Abs[2]{
  \Size[0]{#1}{\lvert}{#2}{}{}{\rvert}
}

\makeatletter

\newcommand{\IfUpperCase}[1]{\begingroup 
  \protected@edef\@tempa{\expandafter\@firstofone\@firstofone#1.}%
  \expandafter\IfUpperCasE \@tempa\delimiter}

\def\IfUpperCasE #1#2\delimiter{%
  \protected@edef\@tempa{\meaning#1\meaning a}%
  \ifnum \expandafter\IfUppercaSE\@tempa \IfUppercaSE
   \endgroup \expandafter\@firstoftwo
  \else
   \endgroup \expandafter\@secondoftwo
  \fi}

\@ifundefined{strip@prefix}{\def\strip@prefix#1>{}}{}
\def\@tempa{the letter }
\edef\@tempa{\expandafter\strip@prefix\meaning\@tempa}

\expandafter\def\expandafter\IfUppercaSE\expandafter#\expandafter1\@tempa#2#3\IfUppercaSE{\uccode`#2=`#2 }

\newif\ifuc@se
\def\setuc@se#1{\IfUpperCase{#1}{\uc@setrue}{\uc@sefalse}}

\def\theoremn@me#1{\ifuc@se \lowercase{\csname#1name\endcsname}\ignorespaces%
  \else \edef\@temp{\lowercase{\lowercase{\csname#1name\endcsname}}}\@temp\ignorespaces%
  \fi}
\def\theoremn@mes#1{\ifuc@se \lowercase{\csname#1names\endcsname}\ignorespaces%
  \else \edef\@temp{\lowercase{\lowercase{\csname#1names\endcsname}}}\@temp\ignorespaces%
  \fi}
  
\def\thmref#1#2{\setuc@se{#1}\lowercase{{\theoremn@me{#1}\lowercase{\ref{#1:#2}}}}}

\newcommand{\DefTheorem}[2]{\newenvironment{#1}[1][\empty]{\ignorespaces\begin{#2}\ifx##1\empty{}\else\lowercase{\label{#1:##1}}\fi\ignorespaces}{\end{#2}\ignorespacesafterend}}

\DefTheorem{Th}{Theorem}
\DefTheorem{Prop}{Proposition}
\DefTheorem{Cor}{Corollary}
\DefTheorem{Lem}{Lemma}
\DefTheorem{Def}{Definition}
\DefTheorem{Rem}{Remark}
\DefTheorem{Par}{para}

\newenvironment{Par*}{\ignorespaces\noindent\ignorespaces}{\ignorespacesafterend}
\makeatletter

\newif\if@smallmat
\newif\if@none
\newif\if@paren
\newif\if@brack
\newif\if@brace
\newif\if@vline
                                 {\ifx#20%
                                        \@smallmattrue%
                                  \else%
                                         \@smallmatfalse
                                  \fi%
                                  \ifx#11%
                                         \@nonefalse\@parentrue\@brackfalse\@bracefalse\@vlinefalse%
                                  \else%
                                       \ifx#12%
                                            \@nonefalse\@parenfalse\@bracktrue\@bracefalse\@vlinefalse%
                                        \else%
                                            \ifx#13%
                                                 \@nonefalse\@parenfalse\@brackfalse\@bracetrue\@vlinefalse%
                                            \else%
                                                 \ifx#14%
                                                       \@nonefalse\@parenfalse\@brackfalse\@bracefalse\@vlinetrue
                                                 \else%
                                                       \ifx#15%
                                                             \@nonefalse\@parenfalse\@brackfalse\@bracefalse\@vlinefalse%
                                                       \else%
                                                             \@nonetrue\@parenfalse\@brackfalse\@bracefalse\@vlinefalse%
                                                       \fi%
                                                 \fi%
                                            \fi%
                                        \fi%
                                   \fi%
                                   \if@smallmat%
                                        \if@none%
                                             \begin{smallmatrix}%
                                        \else%
                                            \if@paren%
                                                  \bigl(\begin{smallmatrix}%
                                            \else%
                                                  \if@brack%
                                                          \bigl[\begin{smallmatrix}%
                                                  \else%
                                                          \if@brace%
                                                               \bigl\{\begin{smallmatrix}%
                                                          \else%
                                                               \if@vline%
                                                                    \bigl\lvert\begin{smallmatrix}%
                                                                \else%
                                                                    \bigl\lVert\begin{smallmatrix}%
                                                                \fi%
                                                          \fi%
                                                  \fi%
                                            \fi%
                                        \fi%
                                   \else%
                                        \if@none%
                                             \begin{matrix}%
                                        \else%
                                            \if@paren%
                                                  \begin{pmatrix}%
                                            \else%
                                                  \if@brack%
                                                          \begin{bmatrix}%
                                                  \else%
                                                          \if@brace%
                                                               \begin{Bmatrix}%
                                                          \else%
                                                               \if@vline%
                                                                    \begin{vmatrix}%
                                                                \else%
                                                                    \begin{Vmatrix}%
                                                                \fi%
                                                          \fi%
                                                  \fi%
                                            \fi%
                                        \fi%
                                   \fi}%
                                  {\if@smallmat%
                                        \if@none%
                                             \end{smallmatrix}%
                                        \else%
                                            \if@paren%
                                                  \end{smallmatrix}\bigr)%
                                            \else%
                                                  \if@brack%
                                                          \end{smallmatrix}\bigr]%
                                                  \else%
                                                          \if@brace%
                                                               \end{smallmatrix}\bigr\}%
                                                          \else%
                                                               \if@vline%
                                                                    \end{smallmatrix}\bigr\rvert%
                                                                \else%
                                                                    \end{smallmatrix}\bigr\rVert%
                                                                \fi%
                                                          \fi%
                                                  \fi%
                                            \fi%
                                         \fi%
                                   \else%
                                        \if@none%
                                             \end{matrix}%
                                        \else%
                                            \if@paren%
                                                  \end{pmatrix}%
                                            \else%
                                                  \if@brack%
                                                          \end{bmatrix}%
                                                  \else%
                                                          \if@brace%
                                                               \end{Bmatrix}%
                                                          \else%
                                                               \if@vline%
                                                                    \end{vmatrix}%
                                                                \else%
                                                                    \end{Vmatrix}%
                                                                \fi%
                                                          \fi%
                                                  \fi%
                                            \fi%
                                        \fi%
                                   \fi}%

\makeatother

\title{Invariant Berezin integration\\ on homogeneous supermanifolds}                                     \author{Alexander Alldridge and Joachim Hilgert\thanks{This research was partly supported by the IRTG ``Geometry and Analysis of Symmetries'' at Universit\"at Paderborn, funded by Deutsche Forschungsgemeinschaft (DFG), Minist\`ere de l'\'Education Nationale (MENESR), and Deutsch-Franz\"osische Hochschule (DFH-UFA).}}
\lastname{Alldridge and Hilgert}


\address{%
Alexander Alldridge\\
Institut f\"ur Mathematik\\
Universit\"at Paderborn\\
Warburger Str.~100\\
33100 Paderborn\\
Germany\\
\texttt{alldridg@math.upb.de}
}

\address{%
Joachim Hilgert\\
Institut f\"ur Mathematik\\
Universit\"at Paderborn\\
Warburger Str.~100\\
33100 Paderborn\\
Germany\\
\texttt{hilgert@math.upb.de}
}

\begin{document}

\maketitle

\begin{abstract}
	Let $\mathcal G$ be a Lie supergroup and $\mathcal H$ a closed subsupergroup. We study the unimodularity of the homogeneous supermanifold $\mathcal G/\mathcal H$, \emph{i.e.}~the existence of $\mathcal G$-invariant sections of its Berezinian line bundle. To that end, we express this line bundle as a $\mathcal G$-equivariant associated bundle of the principal $\mathcal H$-bundle $\mathcal G\to\mathcal G/\mathcal H$. We also study the fibre integration of Berezinians on oriented fibre bundles. As an application, we prove a formula of `Fubini' type: $\int_{\mathcal G}f=(-1)^{\dim\ger h_1\cdot\dim\ger g/\ger h}\int_{\mathcal G/\mathcal H}\int_{\mathcal H}f$, for all $f\in\Gamma_c(G,\mathcal O_{\mathcal G})$. 
	
	Moreover, we derive analogues of integral formulae for the transformation under local isomorphisms $\mathcal G/\mathcal H\to\mathcal S/\mathcal T\!$, and under the products of Lie subsupergroups $\mathcal M\cdot\mathcal H\subset\mathcal U$. The classical counterparts of these formulae have numerous applications in harmonic analysis. 
\end{abstract}

\section{Introduction}

Let $G$ be a Lie group and $H$ a closed subgroup. The homogeneous space $G/H$ is called \emph{unimodular} if there exists a non-zero $G$-invariant volume form. It is an important prerequisite to the study of harmonic and global analysis on the space $G/H$ to find necessary and sufficient conditions for its unimodularity. 

Classically, it is well-known that the unimodularity of $G/H$ is equivalent to the condition that $\det\Ad_{\ger g}=\det\Ad_{\ger h}$ on $H\!$. Put differently, $G/H$ is unimodular if and only if $\bigwedge^{top}(\ger g/\ger h)$ is a trivial $H$-module, where $\ger g$ and $\ger h$ are the Lie algebras of $G$ and $H\!$, respectively.

When $G/H$ is unimodular, one has, for a suitable normalisation of invariant volume forms, the `Fubini' type formula 
\begin{equation}\label{eq:classfubfmla}
	\int_Gf(g)\,dg=\int_{G/H}\int_Hf(gh)\,dh\,d\dot g\mathfa f\in\Cc0G\ .
\end{equation}

Moreover, there is a well-known lemma which allows one to compute the behaviour of invariant-volume forms under local diffeomorphisms of homogeneous spaces for different groups. As an important application, if the Lie group $U$ is as a manifold the direct product of two subgroups $M$ and $H\!$, then for a suitable normalisation of measures, 
\begin{equation}\label{eq:classprodsubgrp}
	\int_U f(u)\,du=\int_{M\times H}f(mh)\frac{\det\Ad_{\ger h}(h)}{\det\Ad_{\ger u}(h)}\,dm\,dh\mathfa f\in\Cc0U\ .
\end{equation}

This formula has manifold applications: It applies, \emph{e.g.}, in the context of Riemannian symmetric spaces, to the Bruhat and Iwasawa decompositions. It plays a role in the proof of the Harish-Chandra isomorphism for Riemannian symmetric spaces, and has many applications in the representation theory of semi-simple Lie groups. 

\medskip\noindent
From the point of view of the geometry and analysis of Lie supergroups and their homogeneous superspaces, it is desirable to have generalisations of all of these facts to the supermathematical context. The basic problem is to characterise unimodularity; here, the Berezinian bundle takes the place of the determinant bundle (whose sections are the volume forms).

We prove that for any Lie supergoup $\mathcal G$ and any closed subsupergroup $\mathcal H$, the Berezinian bundle $\Ber0{\mathcal G/\mathcal H}$ is, as a $\mathcal G$-equivariant vector bundle, isomorphic to the associated bundle $\mathcal G\times^{\mathcal H}\Ber0{(\ger g/\ger h)^*}$ where $\ger g$ and $\ger h$ are, respectively, the Lie superalgebras of $\mathcal G$ and $\mathcal H$ (\thmref{Cor}{quotber}). From this fact, our main theorem (\thmref{Th}{unimodequiv}) follows: $\mathcal G/\mathcal H$ supports a non-zero $\mathcal G$-invariant Berezinian form if and only if $\Ber0{(\ger g/\ger h)^*}$ is a trivial $\mathcal H$-module. Along the way, we discuss all the basic machinery of associated bundles: principal bundles, quotients by free and proper actions, equivariant vector bundles. 

The formula \eqref{eq:classfubfmla} is best understood in the context of fibre integration. We introduce a general fibre integration map for oriented fibre bundles, and show that it satisfies a `Fubini' type formula (\thmref{Prop}{fibreint}). We then apply this to homogeneous principal bundles $\mathcal G\to\mathcal G/\mathcal H$ (\thmref{Prop}{fubformula} and \thmref{Cor}{fubformula}). Finally, in \thmref{Prop}{prodsubsupergrp}, we generalise the formula \eqref{eq:classprodsubgrp}.

\medskip\noindent
Our investigation of invariant Berezin integration is motivated by an ongoing joint project with M.R.~Zirnbauer (K\"oln), concerning the harmonic analysis on a certain class of reductive symmetric superspaces. 

Zirnbauer \cite{zirnbauer-poincare} has employed harmonic superanalysis on the super Poin\-car\'e disk in an application to a problem in mesoscopic physics (the determination of the mean conductivity for a quasi-one dimensional metallic system). These methods fit into a general framework of Riemannian symmetric spaces embedded as subspaces into infinite series of complex symmetric superspaces \cite{zirnbauer-rsss}. In this generality, the harmonic analysis on symmetric superspaces has as yet not been developed. 

In \cite{ahz-chevalley}, we have established a generalisation of Chevalley's restriction theorem to the context of reductive symmetric superpairs. In a series of forthcoming papers, we will employ the results on invariant Berezin integration established in this paper to a generalisation of the Harish-Chandra isomorphism, and to the study of spherical superfunctions on symmetric superspaces. 

\medskip\noindent
Let us fix our notation for what follows. The graded parts of a super vector space $V$ will be written $V_0$ and $V_1$, respectively; $\Pi$ will denote the grading inverting functor. Supermanifolds will in general by denoted $\mathcal X=(X,\mathcal O_{\mathcal X})$, $\mathcal Y=(Y,\mathcal O_{\mathcal Y})$. We will assume, as is common, that the underlying manifolds of supermanifolds are Hausdorff and second countable (where the latter assumption will be used in Section \ref{sec:berezin}). Morphisms of supermanifolds will be denoted by $\vphi=(f,f^*)$ where $f:X\to Y$ and $f^*:\mathcal O_{\mathcal Y}\to f_*\mathcal O_{\mathcal X}$. When appropriate, we will write $\vphi:X\to Y$ and $\vphi^*:\mathcal O_{\mathcal X}\to\vphi_*\mathcal O_{\mathcal Y}$, thereby slightly abusing the notation. We will say that a morphism of supermanifolds is injective, surjective, bijective, open or closed if so is the underlying map of topological spaces. Sometimes we will write $h\in\mathcal O_{\mathcal X}$; by this notation we mean that $h\in\mathcal O_{\mathcal X}(U)$ for some fixed but unspecified open subset $U\subset X$. Finally, for supermanifolds $\mathcal X$ and $\mathcal Y$, let $\mathcal X(\mathcal Y)=\Hom0{\mathcal Y,\mathcal X}$ be the set of morphisms $\mathcal Y\to\mathcal X$. This is also called the set of $\mathcal Y$-points of $\mathcal X$ (by the usual notion that a point is the same as a morphism $*\to\mathcal X$). 

\section{Quotients and actions}

\begin{Par*}
	In this section, we discuss a generalisation of Godement's theorem on quotient manifolds to the context of supermanifolds, due to Almorox \cite{almorox}. As an application, we show that the quotient of a supermanifold by a free and proper Lie supergroup action is again a supermanifold. 
\end{Par*}

\subsection{Quotient supermanifolds}

\begin{Def}
	Consider a morphism of supermanifolds $\vphi:\mathcal X\to\mathcal Y$ given by the map $f:X\to Y$ and the sheaf morphism $f^*:\mathcal O_{\mathcal Y}\to f_*\mathcal O_{\mathcal X}$. We say that $\vphi$ is an \emph{open (closed) embedding} if $f$ is an open (closed) embedding, and $f^*$ is an isomorphism (epimorphism); and that it is a \emph{subsupermanifold} if it factors as the composition of a closed and an open embedding. When the morphism $\vphi$ is understood, then in the latter case, we will sometimes also refer to its domain $\mathcal X$ as a subsupermanifold of $\mathcal Y$. 
	
	Recall $T_x\mathcal X=\Der0{\mathcal O_{\mathcal X,x},\reals}$. The morphism $\vphi$ induces tangent maps $T_x\vphi:T_x\mathcal X\to T_{f(x)}\mathcal Y$ as follows. Given $v\in T_x\mathcal X$ and $h\in\mathcal O_{\mathcal Y,f(x)}$, we have a germ $f^*h\in\mathcal O_{\mathcal X,x}$, and set $[T_x\vphi(v)]h=v(f^*h)$. If $T_x\vphi$ is surjective (injective) \fa $x\in X$, then $\vphi$ is called a \emph{submersion} (an \emph{immersion}). 

	More generally, if $\psi=(g,g^*):\mathcal Z\to\mathcal Y$ is another morphism, then $\vphi$ and $\psi$ are \emph{transversal} whenever for any $x\in X$, $y\in Y$ and $z\in Z$ \scth $f(x)=g(z)=y$, one has $T_y\mathcal Y=T_x\vphi(T_x\mathcal X)+T_z\psi(T_z\mathcal Z)$. If $\vphi$, $\psi$ are transversal, then the fibre product $\mathcal X\times_{\mathcal Y}\mathcal Z$ exists in the category of supermanifolds, and is a submanifold of $\mathcal X\times\mathcal Z$ \cite[Proposition 2.9]{bartocci-bruzzo-hernandez-pestov}.
\end{Def}

\begin{Rem}[superepi]	
	Concerning the definition of submersions, we point out that $\vphi$ being a submersion implies that $\vphi^*:\mathcal O_{\mathcal Y}\to\vphi_*\mathcal O_{\mathcal X}$ is a monomorphism of sheaves \cite[Proposition 2.16.2]{kostant-supergeom}. The converse it obvious, since the graded dimension of $\ger m_x/\ger m_x^2$ equals the graded dimension of $T_x\mathcal X$, for any $x\in X$. Here, $\ger m_x$ denotes the maximal ideal of $\mathcal O_{\mathcal X,x}$. Thus, submersions might have been defined in terms of the structure sheaves.  
\end{Rem}

\begin{Def}
	Let $\mathcal X$ be a supermanifold and $\iota:\mathcal R\to\mathcal X\times\mathcal X$ a subsupermanifold. Let $\delta:\mathcal X\to\mathcal X\times\mathcal X$ be the diagonal morphism and $p_j:\mathcal R\to\mathcal X$ be the projections induced by $\iota$. We call $\mathcal R$ an \emph{equivalence relation} if:
	\begin{enumerate}
		\item There is given a morphism $\vrho:\mathcal X\to\mathcal R$ \scth $\iota\circ\vrho=\delta$,
		\item There is given a morphism $\tau:\mathcal R\times_{\mathcal X}\mathcal R\to\mathcal R$ \scth 
		\[
			p_i\circ\tau=p_i\circ\pi_i\nd\tau\circ(\tau\times\id)=\tau\circ(\id\times\tau)
		\]
		where $\pi_1,\pi_2:\mathcal R\times_{\mathcal X}\mathcal R\to\mathcal R$ are the projections onto the first and second factor, respectively. 
	\end{enumerate}
	
	The fibre product $\mathcal R\times_{\mathcal X}\mathcal R$ exists for the following reason: We have the identity $\iota=(p_1,p_2)$, so that $\rk T_rp_1+\rk T_rp_2=\dim_r\mathcal R$ \fa $r\in R$. By (i), $\vrho$ is an immersion, so $\dim_x\mathcal X\sle\dim_{\vrho(x)}\mathcal R$ for $x\in X$. If $p_1(r)=p_2(r)$, then $r=\vrho(x)$ \fs $x\in X$, and then $p_1$, $p_2$ are transversal at $r$. 
	
	In (ii), the equation for $\tau$ is to be understood on $(\mathcal R\times_{\mathcal X}\mathcal R)\times_{\mathcal X}\mathcal R\cong\mathcal R\times_{\mathcal X}(\mathcal R\times_{\mathcal X}\mathcal R)$. The latter canonical isomorphism obtains since products in any category are commutative, and the fibre product is the product in the category of supermanifolds over $\mathcal X$. 
\end{Def}

\begin{Rem}
	If we were to drop the assumption that $\mathcal R$ be a subsupermanifold of $\mathcal X\times\mathcal X$, then the above axioms would be those of a category with space of objects $\mathcal X$. This will be useful in understanding the treatment, below, of the equivalence relations defined by supergroup actions---these are reminiscent of the `action groupoids' known from the theory of Lie groupoids. 
	
	If $\mathcal R$ is an equivalence relation, then there exists a morphism $\sigma:\mathcal R\to\mathcal R$ \scth $p_1\circ\sigma=p_2$ and $p_2\circ\sigma=p_1$ \cite[Lemma 2.2]{almorox}. 
\end{Rem}

\begin{Def}
	Let $\mathcal X$ be a supermanifold and $\iota:\mathcal R\to\mathcal X\times\mathcal X$ an equivalence relation. If $\vphi:\mathcal X\to\mathcal Y$ is a morphism of supermanifolds, then $\vphi$ is called a \emph{quotient} by $\mathcal R$ if $\vphi$ is a submersion and the coequaliser of $\xymatrix@1@M+1pt{{p_1,p_2\colon\mathcal R}\ar@<2pt>[r]\ar@<-2pt>[r]&{\mathcal X}}$. If moreover, $p_1,p_2$ is the kernel pair of $\vphi$ (\emph{i.e.}~$\iota$ induces an isomorphism $\mathcal R\to\mathcal X\times_{\mathcal Y}\mathcal X$), then $\mathcal R$ is called \emph{effective} and $\mathcal Y$ an \emph{effective quotient}. 
	
	If a quotient of $\mathcal X$ by $\mathcal R$ exists, then it is unique up to canonical isomorphism. In this situation, we write $\mathcal Y=\mathcal X/\mathcal R$. Note that since coequalisers are epimorphisms and the Yoneda embedding preserves colimits, a quotient morphism $\vphi$ is necessarily surjective on the base spaces (consider $\Hom0{*,\vphi}$). 
\end{Def}

\begin{Par*}
	We have the following generalisation to supermanifolds of a theorem which in the context of smooth manifolds is attributed to Godement. 
\end{Par*}

\begin{Th}[quotexist]
	Let $\mathcal X$ be a supermanifold and $\mathcal R$ an equivalence relation on $\mathcal X$. Then $\mathcal X$ admits a quotient by $\mathcal R$ if and only if $\mathcal R$ is closed as a subsupermanifold of $\mathcal X\times\mathcal X$, and the projections $p_i:\mathcal R\to\mathcal X$ are submersions. Whenever a quotient exists, it is effective. 
\end{Th}

\begin{Proof}
	According to \cite[Theorem 2.6]{almorox}, the condition stated in the theorem is necessary and sufficient for the existence of a supermanifold $\mathcal Y$ and a submersion $\vphi:\mathcal X\to\mathcal Y$ \scth the underlying map $\vphi:X\to Y$ is the canonical projection with respect to the equivalence relation $R$ on $X$ underlying $\mathcal R$. 
	
	Assume the condition is satisfied, \emph{i.e.}~that $\mathcal R$ is closed and the $p_i$ are submersions. Then there exists a submersion $\vphi:\mathcal X\to\mathcal Y=(Y,\mathcal O_{\mathcal Y})$ where $Y=X/R$ and $\vphi:X\to Y$ is the canonical projection. Then $\vphi$ is as a morphism of topological spaces the coequaliser of $\xymatrix@1{{p_1,p_2\colon R}\ar@<2pt>[r]\ar@<-2pt>[r]&{X}}$, and the latter are also a kernel pair of $\vphi$ on the level of spaces. By \cite[proof of Theorem 2.6]{almorox}, 
	\[
	\mathcal O_{\mathcal Y}(U)=\Set1{f\in\mathcal O_{\mathcal X}(\vphi^{-1}(U))}{p_1^*f=p_2^*f}
	\]
	and $\vphi^*$ is defined by $\vphi^*f=f$. By this definition, $\vphi^*\colon\mathcal O_{\mathcal Y}\to\vphi_*\mathcal O_{\mathcal X}$ is as a sheaf morphism the equaliser of $\xymatrix@1@M+1pt{{p_1^*,p_2^*\colon\mathcal O_{\mathcal X}}\ar@<2pt>[r]\ar@<-2pt>[r]&{\mathcal O_{\mathcal R}}}$, and the latter form a cokernel pair for $\vphi^*$. Hence, $\vphi$ is an effective quotient morphism. The converse implication is obvious by Almorox's theorem.
\end{Proof}

\begin{Par*}
	Although we shall not use this fact in the sequel, any equivalence relation on a supermanifold in the sense defined above automatically satisfies the assumptions of the theorem in the \emph{odd} variables. To state this more precisely, we make the following definition.
\end{Par*}

\begin{Def}
	A morphism $\vphi=(f,f^*):\mathcal X\to\mathcal Y$ is called an \emph{even (odd) submersion} if \fa $x\in X$, the tangent map $T_x\vphi$ induces by restriction a surjection $(T_x\mathcal X)_i\to(T_{f(x)}\mathcal Y)_i$ where $i=0,1$, respectively. 
\end{Def}

\begin{proposition}[\protect{\cite[Theorem 4.3]{bartocci-bruzzo-hernandez-pestov}}]
	For any equivalence relation $\mathcal R$ on a supermanifold $\mathcal X$, the projections $p_i:\mathcal R\to\mathcal X$ are odd submersions. 
\end{proposition}

\begin{Cor}
	A supermanifold $\mathcal X$ admits a quotient by a given equivalence relation $\mathcal R$ if and only if $\mathcal R$ is closed as a subsupermanifold of $\mathcal X\times\mathcal X$, and the projections $p_i:\mathcal R\to\mathcal X$ are even submersions. 
\end{Cor}

\subsection{Actions of Lie supergroups}

\begin{Def}
	Let $\mathcal G$ be a Lie supergroup and $\mathcal X$ a left $\mathcal G$-space. Let $\alpha:\mathcal G\times\mathcal X\to\mathcal X$ denote the action. Then $\alpha$ is \emph{free} if $(\alpha,\pr_2)\colon\mathcal G\times\mathcal X\to\mathcal X\times\mathcal X$ is an embedding, \emph{transitive} if $(\alpha,\pr_2)$ is a surjective submersion, and \emph{proper} if so is the morphism $(\alpha,\pr_2)$. Here, a morphism of supermanifolds is called \emph{proper} if the underlying map of topological spaces is proper, \emph{i.e.}~closed and with quasi-compact fibres. These definitions can be easily modified for the case of right actions. 
\end{Def}

\begin{Def}
	Let $\mathcal X$, $\mathcal Y$ be supermanifolds with actions $\alpha:\mathcal G\times\mathcal X\to\mathcal X$ and $\beta:\mathcal G\times\mathcal Y\to\mathcal Y$ of a Lie supergroup $\mathcal G$. A morphism $\vphi:\mathcal X\to\mathcal Y$ is \emph{$\mathcal G$-equivariant} if the following diagram commutes:
	\[
		\xymatrix@C+6ex{%
			\mathcal G\times\mathcal X\ar[r]^-{\id\times\vphi}\ar[d]_-{\alpha}&\mathcal G\times\mathcal Y\ar[d]^-{\beta}\\
			\mathcal X\ar[r]_-{\vphi}&\mathcal Y
		}
	\]
	
	We say that $\mathcal X$ admits a \emph{quotient} by $\mathcal G$ if there exists a supermanifold $\mathcal Y$ and a submersion $\pi:\mathcal X\to\mathcal Y$ which is equivariant for the trivial $\mathcal G$-action on $\mathcal Y$, and \scth the following universal property obtains: For any supermanifold $\mathcal Z$, and for any morphism $\vphi:\mathcal X\to\mathcal Z$ which is equivariant for the trivial $\mathcal G$-action, there exists a unique morphism $\psi:\mathcal Y\to\mathcal Z$ making the following diagram commutative:
		\[
		\xymatrix{%
			\mathcal X\ar[rd]_-{\pi}\ar[rr]^-{\vphi}&&\mathcal Z\\
			&\mathcal Y\ar@{-->}[ur]_-{\psi}}
	\] 
	If it exists, $\pi$ is unique up to canonical isomorphism, and we write $\mathcal Y=\mathcal X/\mathcal G$. 
\end{Def}

\begin{Th}[freepropquot]
	Let $\mathcal X$ be a supermanifold and $\mathcal G$ a Lie supergroup acting freely and properly on $\mathcal X$. Then $\mathcal X$ admits a quotient by $\mathcal G$. 
\end{Th}

\begin{Proof}
	Let $\alpha$ denote the action. Then $\iota=(\alpha,\pr_2):\mathcal R=\mathcal G\times\mathcal X\to\mathcal X\times\mathcal X$ is a proper and hence closed embedding. In other words, $\mathcal R$ (or $\iota$) is a closed subsupermanifold of $\mathcal X\times\mathcal X$. 
	
	We need to see that $\mathcal R$ is an equivalence relation with submersive projections. First, note for the projections $p_i:\mathcal R\to\mathcal X$ that $p_2=\pr_2$ whereas $p_1=\alpha$. The required morphism $\vrho:\mathcal X\cong *\times\mathcal X\to\mathcal R$ is given by $\vrho=\eta\times\id$ where $\eta:*\to\mathcal G$ is the unit. Next, we may define $\sigma:\mathcal R\to\mathcal R$ by $\sigma=(i\circ\pr_1,\alpha)$ where $i:\mathcal G\to\mathcal G$ is inversion. Then 
	\[
		p_1\circ\sigma=\alpha\circ(i\circ\pr_1,\alpha)=\pr_2=p_2\ .
	\]
	Moreover, $\sigma^2=\id$ since $i^2=\id$, so that $p_2\circ\sigma=p_1$, and $\sigma$ is an isomorphism. The morphism $p_2=\pr_2$ is manifestly a submersion. Since $\sigma$ is an isomorphism, so is $p_1$.
	
	Next, we consider the fibred product $\mathcal R\times_{\mathcal X}\mathcal R$. We define a morphism
	\[
	\phi:\mathcal R^{(2)}=\mathcal G\times\mathcal G\times\mathcal X\to\mathcal R\times\mathcal R\mathtxt{by}\phi=\Parens1{\pr_1,(\alpha\circ(\pr_2\times\pr_3)),\pr_2\times\pr_3}\ .
	\]
	With the projections $\pi_i:\mathcal R^{(2)}\to\mathcal R$, defined by $\pi_1=(\pr_1,\alpha\circ(\pr_2\times\pr_3))$ and $\pi_2=\pr_2\times\pr_3$, one checks readily that $\mathcal R^{(2)}$ satisfies the universal property of the fibred product of $p_1,p_2$. Thus, we may write $\mathcal R^{(2)}=\mathcal R\times_{\mathcal X}\mathcal R$. Let $\psi=\pr_1\times\pr_3\times\pr_2:\mathcal R\times\mathcal R=\mathcal G\times\mathcal X\times\mathcal G\times\mathcal X\to\mathcal G^2\times\mathcal X$. Then $\psi\circ\phi=\id$, and it follows that $\phi$ is a closed embedding, \emph{i.e.}~a closed subsupermanifold. 
	
	Now, we may define $\tau:\mathcal R^{(2)}\to\mathcal R$ by $\tau=(m\circ(\pr_1\times\pr_2),\pr_3)$ where the morphism $m:\mathcal G\times\mathcal G\to\mathcal G$ is multiplication. Then 
	\[
		p_1\circ\tau=\alpha\circ(m\circ(\pr_1\times\pr_2),\pr_3)=\alpha\circ(\pr_1,\alpha\circ(\pr_2\times\pr_3))=p_1\circ\pi_1
	\]
	and $p_2\circ\tau=\pr_3=p_2\circ\pi_2$. Thus, $\mathcal R$ is indeed an equivalence relation, and the assumptions of \thmref{Th}{quotexist} are satisfied, so that the quotient supermanifold $\pi:\mathcal X\to\mathcal X/\mathcal R$ exists. 
	
	It remains to check that this quotient supermanifold satisfies the universal property of the quotient by a Lie supergroup. To that end, let $\vphi:\mathcal X\to\mathcal Y$ be a morphism which is $\mathcal G$-equivariant with respect to the trivial $\mathcal G$-action on $\mathcal Y$. Then the following diagram commutes:
	\[
		\xymatrix@C+6ex{%
			\mathcal R\ar[r]^-{p_2}\ar[d]_-{p_1}&\mathcal X\ar[d]^-{\vphi}\\
			\mathcal X\ar[r]_-{\vphi}&\mathcal Y
		}
	\]
	Now, by the universal property of coequalisers, there manifestly exists a morphism $\psi:\mathcal X/\mathcal R\to\mathcal Y$ \scth $\vphi=\psi\circ\pi$. This proves the claim. 
\end{Proof}

\section{Principal and associated bundles}

\subsection{Principal bundles}

\begin{Def}
	Given a morphism $p:\mathcal X\to\mathcal B$, the tuple $\mathcal E=(\mathcal X,\mathcal B,p,\mathcal F)$ is called a \emph{fibre bundle with fibre} $\mathcal F$ if there exist an open cover $(\mathcal U_\alpha)$ of $\mathcal B$ and isomorphisms $\tau_\alpha:\mathcal X|_{\mathcal U_\alpha}\to\mathcal U_{\alpha}\times\mathcal F$ \scth the following diagram commutes
	\[
		\xymatrix{%
			\mathcal X|_{\mathcal U_\alpha}\ar[rd]_-{p}\ar[rr]^-{\tau_\alpha}&&\mathcal U_\alpha\times\mathcal F\ar[ld]^-{\pr_1}\\
			&\mathcal U_\alpha}
	\] 
	Here, if $U_\alpha\subset B$ is the underlying manifold of $\mathcal U_\alpha$, $X|_{\mathcal U_\alpha}$ denotes the open subsupermanifold of $\mathcal X$ with base $p^{-1}(U_\alpha)$. The $\tau_\alpha$ are called \emph{local trivialisations}. We shall call $\mathcal X$ the \emph{total space}, $\mathcal B$ the \emph{base space}, $\mathcal F$ the \emph{fibre} and $p$ the \emph{bundle projection} of $\mathcal E$. 

	If $\mathcal G$ is a Lie supergroup and $\mathcal X$, $\mathcal B$ carry left $\mathcal G$-actions $\alpha_{\mathcal X}$, $\alpha_{\mathcal B}$ (say) \scth $p$ is $\mathcal G$-equivariant, then $\mathcal E$ is called a \emph{$\mathcal G$-equivariant fibre bundle}. 
	
	Given a second fibre bundle $\mathcal E'=(\mathcal X',\mathcal B',p',\mathcal F')$, a morphism of fibre bundles $\mathcal E\to\mathcal E'$ is a pair of morphisms $\vphi:\mathcal X\to\mathcal X'$, $f\colon\mathcal B\to\mathcal B'$ \scth $f\circ p=p'\circ\vphi$. We denote the automorphisms of $\mathcal E$ by $\Aut0{\mathcal E}$ or $\Aut[_{\mathcal F}]0{\mathcal X}$. 
\end{Def}

\begin{Def}
	Let $\mathcal X$ be an arbitrary supermanifold, endowed with a left action $\alpha:\mathcal H\times\mathcal X\to\mathcal X$ of a Lie supergroup $\mathcal H$. For any morphism $\gamma:\mathcal Y\to\mathcal H$, define a morphism $\alpha_\gamma:\mathcal Y\times\mathcal X\to\mathcal Y\times\mathcal X$ by $\alpha_\gamma=(\pr_1,\alpha\circ(\gamma\times\id))$. Then $\alpha_\gamma$ is an automorphism of the trivial bundle $\pr_1:\mathcal Y\times\mathcal X\to\mathcal Y$. The action $\alpha$ is called \emph{effective} if $\gamma\mapsto\alpha_\gamma:\mathcal H(\mathcal Y)=\Hom0{\mathcal Y,\mathcal H}\to\Aut[_\mathcal Y]0{\mathcal Y\times\mathcal X}$ is injective for any supermanifold $\mathcal Y$. 
\end{Def}

\begin{Def}[pointsgroupaction]
	Let $\mathcal E=(\mathcal X,\mathcal B,p,\mathcal F)$ be a fibre bundle and $\alpha$ an effective left action of $\mathcal H$ on $\mathcal F$. For any open $U\subset B$, let $\tau_\mathcal X(U)$ be the set of trivialisations of $\mathcal X$ over the open subsupermanifold $\mathcal U\subset\mathcal B$ corresponding to $U$. If $\tau\in\tau_{\mathcal X}(U)$ and $\gamma\in\mathcal H(\mathcal U)$, then $\gamma.\tau=\alpha_\gamma\circ\tau$ defines an effective left action of the group $\mathcal H(\mathcal U)$ on $\tau_{\mathcal X}(U)$ \cite[6.4]{schmitt-supergeom}. 
	
	An \emph{$\mathcal H$-structure} is a subsheaf of sets $\mathcal A\subset\tau_{\mathcal X}$ \scth for every $x\in B$, there exists an open neighbourhood $U$ \scth $\mathcal A(U)\neq\vvoid$, and $\mathcal H(U)$ acts transitively on $\mathcal A(U)$ whenever this set is non-void. \emph{I.e.}, for any two $\tau,\tau'\in\tau_{\mathcal X}(U)$, there exists a $\gamma\in\mathcal H(\mathcal U)$ (unique by effectiveness) \scth $\tau=\gamma.\tau'$. We also say that $\mathcal H$ is the \emph{structure group} of $\mathcal E$, and that the elements of $\mathcal A(U)$ are \emph{bundle charts}. Occasionally, we will also refer to $\tau^{-1}$, where $\tau\in\mathcal A(U)$, as a bundle chart. 
	
	Assume that $\mathcal E$ is at the same time a fibre bundle with $\mathcal H$-structure $\mathcal A$ and a $\mathcal G$-equivariant fibre bundle. Denote the $\mathcal H$-action in the fibre by $\alpha$, and the $\mathcal G$-actions on $\mathcal X$ and $\mathcal B$ by $\tilde\beta$ and $\beta$, respectively. We say that $\mathcal E$ is a \emph{$\mathcal G$-equivariant fibre bundle with $\mathcal H$-structure} if for each $b\in B$, there are open subsupermanifolds $\mathcal U,\mathcal W_i$ of $\mathcal B$ and $\mathcal V_i$ of $\mathcal G$ \scth $b\in U$, $(\mathcal V_i\times\mathcal W_i)$ is an open cover of $\beta^{-1}(\mathcal U)$, and there exist bundle charts $\tau\in\mathcal A(U)$, $\tau_i\in\mathcal A(W_i)$, and $\gamma_i\in\mathcal H(\mathcal V_i\times\mathcal W_i)$ \scth \fa $i$, the following diagram commutes
	\begin{equation}\label{eq:gequivarhbundledef}
		\xymatrix@C+6ex@R-3ex{%
			\mathcal V_i\times\mathcal X|_{\mathcal W_i}\ar[r]^-{\id\times\tau_i}\ar[dd]_-{\tilde\beta}
			&\mathcal V_i\times\mathcal W_i\times\mathcal F\ar[dr]^-{(\id,\gamma_i)\times\id}\\
			&&\mathcal V_i\times\mathcal W_i\times\mathcal H\times\mathcal F\ar[dl]^-{\beta\times\alpha}\\
			\mathcal X|_{\mathcal U}\ar[r]_-{\tau}&\mathcal U\times\mathcal F
		}
	\end{equation}
	Equivalently, the morphism $\mathcal G\times\mathcal X\to\beta^*\mathcal X$ induced by $\tilde\beta$ is one of $\mathcal H$-fibre bundles, in the sense of \cite[6.6]{schmitt-supergeom}. 
\end{Def}

For later use and reference, we record the following fundamental fact \cite[Propositions 5.3, 6.5]{schmitt-supergeom}.

\begin{Prop}
	Let $\mathcal B$, $\mathcal F$ be supermanifolds, $\mathcal H$ a Lie supergroup acting effectively on $\mathcal F$, $(\mathcal U_i)$ an open cover of $\mathcal B$, and $\mathcal U_{ij}=\mathcal U_i\cap\mathcal U_j$. Assume given $\vphi_{ij}\in\mathcal H(\mathcal U_{ij})$ \scth $\vphi_{ij}\cdot\vphi_{jk}=\vphi_{ik}$ on $\mathcal U_{ijk}=\mathcal U_{ij}\cap\mathcal U_{jk}$ and $\vphi_{ii}=1$. Here, $\vphi_{ij}\cdot\vphi_{jk}=m\circ(\vphi_{ij}\times\vphi_{jk})$ and $1$ denotes the unique morphism $\mathcal U_{ii}=\mathcal U_i\to\mathcal H$ which factors through the unit $*\to\mathcal H$. 
	
	There exists an $\mathcal H$-fibre bundle $\mathcal E=(\mathcal X,\mathcal B,p,\mathcal F)$ with bundle charts $\tau_i$ over $\mathcal U_i$ \scth $\vphi_{ij}.\tau_j=\tau_i$ on $\mathcal U_{ij}$, and $\mathcal E$ is unique up to unique isomorphism. We say that $\mathcal E$ is \emph{determined by the cocyle $(\vphi_{ij})$}. 
\end{Prop}

\begin{Def}
	Let $\mathcal E=(\mathcal X,\mathcal B,p,\mathcal H)$ be a fibre bundle with fibre and structure group $\mathcal H$ where $\mathcal H$ a Lie supergroup which acts from the left on itself via left multiplication. (This action is effective \cite[Lemma 6.3]{schmitt-supergeom}.) Then $\mathcal E$ is called a \emph{principal $\mathcal H$-bundle}. 
\end{Def}

\begin{Prop}[principalbundle]
	Let $\mathcal X$ be a supermanifold and $\mathcal H$ a Lie supergroup. If $\mathcal H$ acts freely and properly from the right on $\mathcal X$ via $\alpha:\mathcal X\times\mathcal H\to\mathcal X$, then $\mathcal E=(\mathcal X,\mathcal B,p,\mathcal H)$ is a principal $\mathcal H$-bundle, where $p:\mathcal X\to\mathcal B$ is the quotient morphism. Conversely, given any morphism $p:\mathcal X\to\mathcal B$ \scth the tuple $\mathcal E=(\mathcal X,\mathcal B,p,\mathcal H)$ is a principal $\mathcal H$-bundle, the supermanifold $\mathcal X$ allows for a free and proper $\mathcal H$-action \scth $p$ is the quotient morphism.
\end{Prop}

\begin{Proof}
	Assume that $\mathcal E$ be a principal $\mathcal H$-bundle. For each local trivialisation $\tau\in\tau_{\mathcal X}(U)$, $\mathcal X|_{\mathcal U}\cong\mathcal U\times\mathcal H$ is endowed with a right $\mathcal H$-action which is induced by the canonical right action of $\mathcal H$ on itself. Since the latter commutes with the canonical left action, we obtain a right action of $\mathcal H$ on $\mathcal X$ \cite[Proposition 6.18]{schmitt-supergeom}. It is clear from its definition that $p$ satisfies the universal property of the quotient morphism when restricted to $\mathcal X|_{\mathcal U}$, and hence globally. 
	
	On the level of ordinary manifolds, $E=(X,B,p,H)$ is a principal $H$-bundle, so the action of $H$ on $X$ is free and proper. To check that the morphism $(\pr_1,\alpha):\mathcal X\times\mathcal H\to\mathcal X\times\mathcal X$ is an embedding, it suffices to prove that it is an immersion. We may assume that $\mathcal X=\mathcal B\times\mathcal H$, since this is a local property. But the multiplication morphism $m:\mathcal H\times\mathcal H\to\mathcal H$ is a free $\mathcal H$-action; indeed, $(\pr_1,m)$ is an isomorphism with inverse $(\pr_1,m\circ(i\times\id))$ where $i$ denotes inversion. This implies that $(\pr_1,\alpha)$ is a proper embedding, so that $\alpha$ is free and proper. 
	
	Conversely, assume that $\alpha$ be a free and proper $\mathcal H$-action. Then the quotient $\mathcal X/\mathcal H$ exists, and we assume that $p:\mathcal X\to\mathcal B$ be the quotient morphism. Since $p$ is a surjective submersion, $\mathcal B$ has an open cover $(\mathcal U_\alpha)$ \scth there are morphisms $s_\alpha:\mathcal U_\alpha\to\mathcal X$ \scth $p\circ s_\alpha=\id_{\mathcal U_\alpha}$. It follows that the $s_\alpha$ are open embeddings. Define $\phi_\alpha:\mathcal U_\alpha\times\mathcal H\to\mathcal X|_{\mathcal U_\alpha}$ by $\phi_\alpha=\alpha\circ(s_\alpha\times\id)$. This is well-defined because $p\circ\alpha=p\circ\pr_1$. 
	
	We have $p\circ\phi_\alpha=p\circ s_\alpha\circ\pr_1=\pr_1$. It is clear that $\phi_\alpha$ is bijective on the level of ordinary spaces, and it is an immersion. Moreover, since $\mathcal B$ is an effective quotient by \thmref{Th}{quotexist}, we have $\mathcal X\times_{\mathcal B}\mathcal X\cong\mathcal R\cong\mathcal X\times\mathcal H$ (\emph{cf.}~\thmref{Th}{freepropquot} and proof) and thus locally $2\dim\mathcal X-\dim\mathcal B=\dim\mathcal X+\dim\mathcal H$; hence, $\dim\mathcal X=\dim\mathcal B+\dim\mathcal H$. Thus, $\phi_\alpha$ is an isomorphism by \cite[Theorem 2.16 and Corollary]{kostant-supergeom}. It follows that the $\phi_\alpha$ define local trivialisations of $\mathcal E$, and they are manifestly $\mathcal H$-equivariant. That the $\phi_\alpha$ define an $\mathcal H$-structure by considering the left action of $\mathcal H$ on itself follows from \cite[Lemma 6.2 and Proposition 6.5]{schmitt-supergeom}. 
\end{Proof}

\subsection{Associated bundles}

\begin{Def}
	Let $\mathcal X,\mathcal Y$ be supermanifolds and $\mathcal H$ a Lie supergroup, and assume given a right action $\alpha:\mathcal X\times\mathcal H\to\mathcal X$ and a left action $\beta:\mathcal H\times\mathcal Y\to\mathcal Y$ where $\alpha$ is free and proper. Consider the diagonal right action $\gamma:\mathcal X\times\mathcal Y\times\mathcal H\to\mathcal X\times\mathcal Y$, defined by $\gamma=(\alpha\circ(\pr_1,\pr_3),\beta\circ(i\circ\pr_3,\pr_2))$. 
	
	It is clear that $(\gamma,\pr_1,\pr_2)$ is a proper injection on the level of ordinary spaces, and since $T_{(x,h)}\alpha$ is injective for any $(x,h)\in X\times H$, it follows that $\gamma$ is an immersion. Hence, $\gamma$ is a free and proper action, and the quotient 
	\[
	\mathcal X\times^{\mathcal H}\mathcal Y=(\mathcal X\times\mathcal Y)/\mathcal H
	\]
	exists. By the universal property, the projection $\pr_1:\mathcal X\times\mathcal Y\to\mathcal X$ descends to a submersion $p:\mathcal X\times^{\mathcal H}\mathcal Y\to\mathcal X/\mathcal H$ which we will call the \emph{induced projection}. 
	
	Assume given a local trivialisation of $\pi:\mathcal X\to\mathcal X/\mathcal H$, $\tau:\mathcal X|_{\mathcal U}\to\mathcal U\times\mathcal H$ (say). Then $(\mathcal U\times\mathcal H)\times^{\mathcal H}\mathcal Y\cong\mathcal U\times\mathcal Y$, and this defines a local trivialisation of $p$. 
\end{Def}

\begin{Prop}[assocbundle]
	Given a principal $\mathcal H$-bundle $\pi:\mathcal X\to\mathcal B$ and a left $\mathcal H$-supermanifold $\mathcal Y$, the tuple $\mathcal E=(\mathcal X\times^{\mathcal H}\mathcal Y,\mathcal B,p,\mathcal Y)$ where $p:\mathcal X\times^{\mathcal H}\mathcal Y\to\mathcal B$ is the induced projection, is a fibre bundle, called the bundle with fibre $\mathcal Y$ \emph{associated with} $\mathcal X$. If the action of $\mathcal H$ on $\mathcal Y$ is effective, then $\mathcal X\times^{\mathcal H}\mathcal Y$ has an $\mathcal H$-structure.
\end{Prop}

\begin{Proof}
	It only remains to specify the $\mathcal H$-structure if the action on $\mathcal Y$ is effective. With $\tau\in\tau_{\mathcal X}(U)$, we associate $\tilde\tau\in\tau_{\mathcal X\times^{\mathcal H} V}(U)$, the morphism induced on the quotient by the composite 
	\begin{equation}\label{eq:assochstructure}
		\xymatrix@C+6ex@M+2pt{%
			\tau'\colon\mathcal X|_{\mathcal U}\times\mathcal Y\ar[r]^-{\tau\times\id}
			&\mathcal U\times\mathcal H\times\mathcal Y\ar[r]^-{\id\times\alpha}
			&\mathcal U\times\mathcal Y
		}
	\end{equation}
	where $\alpha$ is the action of $\mathcal H$ on $\mathcal Y$. This defines a subsheaf of $\tau_{\mathcal X\times^{\mathcal H}\mathcal Y}(U)$. 
	
	Recall the action of $\mathcal H(U)$ on local trivialisations from \thmref{Def}{pointsgroupaction}. Let $\sigma,\tau\in\tau_{\mathcal X}(U)$, there exists $\gamma\in\mathcal H(U)$ \scth $\sigma=\gamma.\tau=m_\gamma\circ\tau$ where $m$ is the multiplication of $\mathcal H$. We compute
	\begin{align*}
		\gamma.\tau'&=(\pr_1,\alpha\circ(\gamma\times\id))\circ(\id\times\alpha)\circ(\tau\times\id)\\
		&=(\id\times\alpha)\circ(\pr_1,m\circ(\gamma\circ\pr_1,\pr_2),\pr_3)\circ(\tau\times\id)=\sigma'\ .
	\end{align*}
	It follows that $\gamma.\tilde\tau=\tilde\sigma$.
\end{Proof}

\begin{Def}
A fibre bundle $\mathcal E$ whose fibre is a \emph{linear} supermanifold $V$ and which is endowed with a $\mathcal{GL}(V)$-structure is called a \emph{vector bundle}. One might also define vector bundles as locally free $\mathcal O_{\mathcal B}$-module sheaves (where $\mathcal B$ is the base space of $\mathcal E$). These notions are equivalent by \cite[Proposition 7.33]{schmitt-supergeom}. 

Let $\mathcal G$ be a Lie supergroup. If $\mathcal E$ is a vector bundle endowed with $\mathcal G$-actions $\tilde\alpha$, $\alpha$ on $\mathcal X$ and $\mathcal B$, respectively, then $\mathcal E$ is called a \emph{$\mathcal G$-equivariant vector bundle} if it is a $\mathcal G$-equivariant fibre bundle with $\mathcal{GL}(V)$-structure in the sense of \thmref{Def}{pointsgroupaction}. 

We shall also need the notion of a morphism of vector bundles. To that end, let $\mathcal E=(\mathcal X,\mathcal B,p,V)$ and $\mathcal E'=(\mathcal X',\mathcal B,p',V')$ be vector bundles over the same base. A \emph{morphism of vector bundles} $\mathcal E\to\mathcal E'$ is given by a bundle morphism $f:\mathcal X\to\mathcal X'$ \scth for each $b\in B$, there are an open neighbourhood $U\subset B$ of $x$, vector bundle charts $\sigma\in\mathcal A_{\mathcal X}(U)$, $\tau\in\mathcal A_{\mathcal X'}(U)$, and a morphism $\vphi:\mathcal U\to\Hom0{V,V'}_0$ \scth the following diagram commutes:
	\begin{equation}\label{eq:vbmordef}
		\xymatrix@C+6ex{%
			\mathcal X|_{\mathcal U}\ar[rr]^-{f}\ar[d]_-{\tau}
			&&\mathcal X'|_{\mathcal U}\ar[d]^-{\sigma}\\
			\mathcal U\times V\ar[rd]_-{(\id,\vphi)\times\id\quad}
			&&\mathcal U\times V'\\
			&\mathcal U\times\Hom0{V,V'}_0\times V\ar[ru]_-{\id\times\eps}
		}
	\end{equation}
	where $\eps:\Hom0{V,V'}_0\times V\to V'$ is the natural (`evaluation') morphism. 
	
	If $\mathcal E'$ is a vector bundle over a different base $\mathcal B'$ (say), then a \emph{morphism of vector bundles} $\mathcal E\to\mathcal E'$ is a pair of morphisms $\vphi:\mathcal X\to\mathcal X'$, $\phi:\mathcal B\to\mathcal B'$ \scth $p'\circ\vphi=\phi\circ p$ and $\vphi$ induces a vector bundle morphism $\mathcal X\to\phi^*\mathcal X'$ (where the pullback is defined in \cite[Propositions 5.6, 7.21]{schmitt-supergeom}).
\end{Def}

\begin{Prop}
	Let $\mathcal X$ be a principal $\mathcal H$-bundle and $V$ be endowed with a linear $\mathcal H$-action. Then the associated bundle with fibre $V$ is a vector bundle. 
\end{Prop}

\begin{Proof}
	The action of $\mathcal H$ on $V$ factors through a morphism of supergroups $\mathcal H\to\mathcal{GL}(V)$. The same proof as that of the existence of $\mathcal H$-structures in \thmref{Prop}{assocbundle} shows that $\mathcal X$ has a $\mathcal{GL}(V)$-structure. (It is easy to check by the definition that the canonical action of $\mathcal{GL}(V)$ on $V$ is effective.)
\end{Proof}

\section{Homogeneous superspaces}

In what follows, let $\mathcal G$ be a Lie supergroup and $\mathcal H$ a closed sub-supergroup. 

\subsection{The tangent bundle of $\mathcal G/\mathcal H$ as an associated bundle}\ \\

\begin{Par}
	We recall some basic facts related to tangent morphisms. Let $\mathcal X$ be a supermanifold. The tangent bundle $T\mathcal X\to\mathcal X$ is the vector bundle which is associated with the locally free $\mathcal O_{\mathcal X}$-module $\Der0{\mathcal O_{\mathcal X}}$. 
	
	Define $\Omega^1_{\mathcal X}$ by $\Pi\Omega^1_{\mathcal X}=\Der0{\mathcal O_{\mathcal X}}^*$. If $f:\mathcal X\to\mathcal Y$ is a morphism, then by \cite[Proposition 8.12]{schmitt-supergeom}, there is a unique sheaf morphism $f^*:\Omega^1_{\mathcal Y}\to f_*\Omega^1_{\mathcal X}$ \scth $f^*(d\vphi)=df^*\vphi$ \fa $\vphi\in\mathcal O_{\mathcal Y}$. Thus, there is a vector bundle morphism $Tf:T\mathcal X\to T\mathcal Y$ where $(Tf)^*\Pi(\omega)=\Pi(f^*\omega)$. 
	
	Put differently, the set of local sections of $T\mathcal X$ is $\Der0{\mathcal O_{\mathcal X}}$. Then we say that $v\in\Der0{\mathcal O_{\mathcal Y}}$ is \emph{$f$-related} to $u\in\Der0{\mathcal O_{\mathcal X}}$ if and only if $f^*\omega(u)=\omega(v)$ \fa $\omega\in\Omega^1_{\mathcal Y}$, if and only if $d(f^*\vphi)(u)=d\vphi(v)$ \fa $\vphi\in\mathcal O_{\mathcal Y}$. Recall from \cite[7.27]{schmitt-supergeom} that we may interpret sections of $T\mathcal X$ as morphisms $\mathcal X\to T\mathcal X\oplus\Pi T\mathcal X$ which are right inverses of the bundle projection $T\mathcal X\oplus\Pi T\mathcal X\to\mathcal X$ (we have to add the tangent bundle with the opposite parity in order to treat even and odd sections). If we do so, then $v$ is $f$-related to $u$ if and only if $v\circ f=(Tf\oplus\Pi Tf)\circ u$. This also determines $Tf$.
	
	The assignment $\mathcal X\mapsto T\mathcal X$, $f\mapsto Tf$ defines a product-preserving functor from supermanifolds to vector bundles \cite[8.16-17]{schmitt-supergeom}, the \emph{tangent functor}. By functoriality, if $\mathcal G$ is a Lie supergroup, so is $T\mathcal G$, and $\mathcal G$ is a closed subsupergroup via the zero section; if $\mathcal X$ is a space with a $\mathcal G$-action, then $T\mathcal X$ has a $T\mathcal G$-action, and in particular, $T\mathcal X$ is a $\mathcal G$-equivariant vector bundle. 
\end{Par}

\begin{Prop}
	The tangent bundle $T\mathcal G$ is $\mathcal G$-equivariantly trivial. More precisely, the composite
	\[
		\xymatrix@C+4ex@M+2pt{%
			\mathcal G\times\ger g\ar[r]^-{0\times\iota}&T\mathcal G\times T\mathcal G\ar[r]^-{Tm}&T\mathcal G
		}
	\]
	is a $\mathcal G$-equivariant vector bundle isomorphism. Here, $0:\mathcal G\to T\mathcal G$ is the zero section, and $\iota:\ger g=T_1\mathcal G\to T\mathcal G$ is the inclusion of the fibre at the unit. 
\end{Prop}

\begin{Proof}
	Let $\phi$ denote the morphism defined in the assertion. The zero section exhibits $\mathcal G$ as a closed subsupergroup of $T\mathcal G$, and by restriction of $Tm$ to from $T\mathcal G^2$ to $\mathcal G\times T\mathcal G$, we obtain the natural left $\mathcal G$-action on $T\mathcal G$. The action of $\mathcal G$ on $\mathcal G\times\ger g$ is simply $m\times\id$. From these definitions, it is clear that $Tm$, and hence $\phi$, is $\mathcal G$-equivariant.
	
	Let $p:T\mathcal G\to\mathcal G$ be the bundle projection. Now, consider the morphism $\psi=(Ti\circ\pr_1,Tm)\circ(Ti\times\id)\circ(0\circ p,\id):T\mathcal G\to T\mathcal G$. We compute
	\begin{align*}
		\psi\circ\phi&=(Ti\circ\pr_1,Tm)\circ(\id\times Tm)\circ(Ti\times\id\times\id)\circ(\delta\times\id)\circ(0\times\iota)\\
		&=(\pr_1,Tm\circ(Tm\times\id)\circ(Ti\times\id\times\id)\circ(\delta\times\id))\circ(0\times\iota)\\
		&=0\times\iota\ .
	\end{align*}
	It follows that $\phi$ is an injective immersion. Moreover, 
	\[
	p\circ\phi=m\circ(p\times p)\circ(0\times\iota)=m\circ(\id\times\eta)=\pr_1
	\]
	where we also write $\eta$ for the unique morphism $\ger g\to\mathcal G$ which factors through $\eta:*\to\mathcal G$. Thus, $\phi$ is a vector bundle morphism along the identity. We conclude that $\phi$ is an isomorphism.
\end{Proof}

\begin{Par}
	Let $\pi:\mathcal G\to\mathcal G/\mathcal H$ be the quotient morphism and denote the unit tangent spaces by $\ger g=T_1\mathcal G$, $\ger h=T_1\mathcal H$. Let $\ker T\pi$ be the kernel (in the category of vector bundles over $\mathcal G$) of the canonical map $T\mathcal G\to\pi^*T(\mathcal G/\mathcal H)$. Any bundle chart $\tau:T\mathcal G|_{\mathcal U}\to\mathcal U\times T_x\mathcal G$ of $T\mathcal G$ restricts to a bundle chart $(\ker T\pi)|_{\mathcal U}\to\mathcal U\times\ker T_x\pi$. 
\end{Par}

\begin{Lem}
	The subbundle $\ker T\pi\subset T\mathcal G$ identifies with $\mathcal G\times\ger h$. Therefore, $T\pi$ induces a $\mathcal G$-equivariant vector bundle morphism 
	\[
	\mathcal G\times\ger g/\ger h=T\mathcal G/\ker T\pi\to T(\mathcal G/\mathcal H)\ .
	\]
\end{Lem}

\begin{Proof}
	Let $\vphi\in\mathcal O_{\mathcal G/\mathcal H}$. Then $f=\pi^*\vphi\in\mathcal O_{\mathcal G}$ is characterised by $m^*f=\pr_1^*f$ for $m,\pr_1:\mathcal G\times\mathcal H\to\mathcal G$. If $v$ is $\phi$-related to a local section $u$ of $G\times\ger h$, then 
	\[
		df(v)=dm^*f((0\times\iota)(u))=d\pr_1^*f((0\times\iota)(u))=0\ .
	\]
	It follows that $v$ is a local section of $\ker T\pi$. By equality of dimension $\phi$ induces a $\mathcal G$-equivariant isomorphism $\mathcal G\times\ger h\cong\ker T\pi$. 
\end{Proof}

\begin{Par}
	On $T\mathcal G$, consider the right action $\alpha$ of $\mathcal H$ induced by the tangent multiplication $Tm:T\mathcal G\times T\mathcal H\to T\mathcal G$. Thus, $\alpha=Tm\circ(\id\times0):T\mathcal G\times\mathcal H\to T\mathcal G$. Recall also the adjoint action $\Ad:\mathcal G\times\ger g\to\ger g$. It is given as the restriction of the composite 
	\[
		\xymatrix@C+3ex{%
			T\mathcal G^2\ar[r]^-{\delta\times\id}&T\mathcal G^3\ar[r]^-{\id\times Ti\times\id}&T\mathcal G^3\ar[r]^-{(23)}&T\mathcal G^3\ar[r]^-{Tm^{(2)}}&T\mathcal G}
	\]
	where $m^{(2)}=m\circ(m\times\id)$ and $(23)$ interchanges the second and third factor.
\end{Par}

\begin{Lem}
	Let $\gamma$ be the diagonal right $\mathcal H$-action on $\mathcal G\times\ger g$. The isomorphism $\phi:\mathcal G\times\ger g\to T\mathcal G$ is $\mathcal H$-equivariant.
\end{Lem}	

\begin{Proof}
	Let $\delta^{(2)}=(\delta\times\id)\circ\delta$ and $\id^{(2)}=\id\times\id$. Let $\tilde\gamma$ be the morphism
	\[
		\xymatrix@C+7ex{%
			T\mathcal G^3\ar[r]^-{\id^{(2)}\times\delta^{(2)}}&T\mathcal G^5\ar[r]^-{(34)\circ(23)}&T\mathcal G^5\ar[r]^-{\id^{(2)}\times Ti\times\id^{(2)}}&T\mathcal G^5\ar[r]^-{Tm\times Tm^{(2)}}&T\mathcal G^2
		}
	\]
	
	Then $\gamma$ satisfies $(0\times\iota)\circ\gamma=\tilde\gamma\circ(0\times\iota\times0)$, and this determines the morphism $\gamma$ uniquely. But $\phi=Tm\circ(0\times\iota)$ and $Tm\circ\gamma'=Tm^{(2)}$. This proves the assertion. 
\end{Proof}

\begin{Prop}[quottangent]
The morphism $\mathcal G\times\ger g/\ger h=T\mathcal G/\ker T\pi\to T(\mathcal G/\mathcal H)$ induced by $T\pi$ is the quotient morphism for the induced right $\mathcal H$-action. In particular, $T(\mathcal G/\mathcal H)$ is $\mathcal G$-equivariantly isomorphic, as a vector bundle, to $\mathcal G\times^{\mathcal H}\ger g/\ger h$. 
\end{Prop}

\begin{Proof}
	For the trivial $\mathcal H$-action on $\mathcal G/\mathcal H$, $\pi$ is $\mathcal H$-equivariant. Applying the tangent functor, $T\pi$ is $T\mathcal H$-equivariant for the trivial $T\mathcal H$-action on $T(\mathcal G/\mathcal H)$. In particular, $T\pi$ is $\mathcal H$-equivariant, and hence, so is the morphism induced by $T\pi$, namely $\widetilde{T\pi}:T\mathcal G/\ker T\pi\to T(\mathcal G/\mathcal H)$.
	
	Hence, there exists a vector bundle morphism $\vphi:\mathcal G\times^{\mathcal H}\ger g/\ger h\to T(\mathcal G/\mathcal H)$ \scth $\vphi\circ p=\widetilde{T\pi}\circ\tilde\phi$ where $\tilde\phi:\mathcal G\times\ger g/\ger h\to T\mathcal G/\ker T\pi$ is the isomorphism induced by $\phi$, and $p:\mathcal G\times\ger g/\ger h\to\mathcal G\times^{\mathcal H}\ger g/\ger h$ is the quotient morphism. We know that $p$ and $\widetilde{T\pi}$ are along $\pi$, and $\tilde\phi$ is along $\id$. On the other hand, $\vphi$ is a submersion since so is $\widetilde{T\pi}\circ\tilde\phi$. Hence, $\vphi$ is an isomorphism, and by definition, it is $\mathcal G$-equivariant. 
\end{Proof}

\subsection{Berezinians}\ \\

\begin{Par}
	Let $V$ be a finite-dimensional super-vector space. Recall that $\mathcal{GL}(V)$ is the open subsupermanifold of the linear supermanifold $\End0V_0$ corresponding to the open subset $\GL(V)\subset\End0V_0$. Moreover, for any supermanifold $\mathcal Y$, the $\mathcal Y$-points of $\mathcal{GL}(V)$ are $\Aut[_{\mathcal O(\mathcal Y)}]0{\mathcal O(\mathcal Y)\otimes V}$, and this induces the Lie supergroup structure on $\mathcal{GL}(V)$.
	
	The standard action of $\mathcal{GL}(V)$ on $V$ is defined by the actions 
	\[
		\alpha(\mathcal Y):\Aut[_{\mathcal O(\mathcal Y)}]0{\mathcal O(\mathcal Y)\otimesÊV}\times(\mathcal O(\mathcal Y)\otimes V)_0\to(\mathcal O(\mathcal Y)\otimes V)_0
	\] 
	which are natural in $\mathcal Y$. One defines 
	\[
	\alpha(\mathcal Y)^*:\Aut[_{\mathcal O(\mathcal Y)}]0{\mathcal O(\mathcal Y)\otimesÊV}\times(\mathcal O(\mathcal Y)\otimes V^*)_0\to(\mathcal O(\mathcal Y)\otimes V^*)_0
	\]
	by 
	\[
		\alpha(\mathcal Y)^*(\gamma,f\otimes u^*)(g\otimes v)=(f\otimes u^*)\Parens1{\alpha(\mathcal Y)(\gamma^{-1},g\otimes v)}\ .
	\]
	
	This is again natural in $\mathcal Y$, and defines the contragredient action of $\mathcal{GL}(V)$ on $V^*$. Moreover, we recall that there is a natural Lie supergroup isomorphism $\mathcal{GL}(V)\to\mathcal{GL}(\Pi V)$ given on the level of $\mathcal Y$-points by $\gamma\mapsto\Pi\gamma\Pi$. In particular, the identity $\Pi:V\to\Pi V$ is $\mathcal{GL}(V)$-equivariant. Using $\mathcal Y$-points, one shows equally easily that the canonical isomorphism $\Pi V\otimes V^*\cong\Hom0{V,\Pi V}$ is $\mathcal{GL}(V)$-equivariant.
\end{Par}\ \\

\begin{Par}[berpbdef]
	Let $V$ be a finite-dimensional super-vector space, and consider the (super-) symmetric algebra $S(\Pi V\oplus V^*)$. The identity $\Pi:V\to\Pi V$ may be considered as an element of $\Pi V\otimes V^*$ and thus embedded as an odd element of $S(\Pi V\oplus V^*)$. Since $\Pi$ is odd, $\Pi\cdot\Pi=0$ in $S(\Pi V\oplus V^*)$, and multiplication by $\Pi$ is a differential on this vector space. One defines $\Pi^p(\Ber0V)$, where $p=\dim V_0$, to be the homology of this differential. One can show that $\Ber0V$ is a super-vector space of total dimension one, and parity $q=\dim V_1$. Thus, $\dim\Ber0V=1|0$ if $q$ is even, and $\Ber0V=0|1$ if $q$ is odd. Since $\Pi$ is $\mathcal{GL}(V)$-equivariant, $\Ber0V$ carries a linear $\mathcal{GL}(V)$-action. We shall call $\Ber0V$ the \emph{Berezinian module} and the corresponding action the \emph{Berezinian action}. In fact, this definition can be analogously performed whenever $V$ is a graded free and finitely generated $R$-module, where $R$ is any commutative ring. We will use this fact in one instance below, and to stress the base ring, we will then write $\Ber[_R]0V$. 
	
	Let $\mathcal E=(\mathcal X,\mathcal B,p,V)$ be a vector bundle. We define the \emph{Berezinian bundle} as follows. Let $\mathcal A$ be the $\mathcal{GL}(V)$-structure of $\mathcal E$ and $\tau_i\in\mathcal A(U_i)$ bundle charts over some open cover $(U_i)$ of $B$. Let $\vphi_{ij}:\mathcal U_{ij}\to\mathcal{GL}(V)$ be the corresponding cocyle \scth $\tau_i=\vphi_{ij}.\tau_j$. By application of the Berezinian, we define a cocyle $\psi_{ij}=\Ber0{\vphi_{ij}}:\mathcal U_{ij}\to\mathcal{GL}(\eps)$, where $\eps=1|0$ or $\eps=0|1$ according to the parity of $q$. The corresponding vector bundle on $\mathcal B$ with fibre $\Pi^q(\reals)$ is the Berezinian $\Ber0{\mathcal E}=\Ber[_{\mathcal B}]0{\mathcal X}$, cf.~\cite[Propositions 5.3, 6.5]{schmitt-supergeom}. (The definition has to be modified appropriately when $\mathcal X$ does not have pure dimension.) 
	
	If $\mathcal X$ is any supermanifold, then we write $\Ber0{\mathcal X}=\Ber0{T^*\mathcal X}$. Here, $T^*\mathcal X$ is the dual bundle of $T\mathcal X$. Whenever $\phi:\mathcal X\to\mathcal Y$ is an isomorphism of supermanifolds, we denote $T^*\phi:T^*\mathcal Y\to T^*\mathcal X$ the super-transpose of $T\phi$. Hence, there is an induced isomorphism $\Pi(T^*(\phi^{-1}))\oplus T\phi:\Pi(T^*\mathcal X)\oplus T\mathcal X\to\Pi(T^*\mathcal Y)\oplus T\mathcal Y$ of vector bundles. This induces an isomorphism $\Ber0{\mathcal X}\to\Ber0{\mathcal Y}$ of line bundles which we denote by $\Ber0\phi$. For any local section $\omega\in\Gamma(U,\Ber0{\mathcal Y})$, we define the \emph{pullback} $\phi^*\omega\in\Gamma(\phi^{-1}(U),\Ber0{\mathcal X})$ by $\Ber0\phi\circ\phi^*\omega=\omega\circ\phi$. (Recall that we may consider a section of a vector bundle $\mathcal E\to\mathcal X$ as a morphism $\mathcal X\to\mathcal E\oplus\Pi(\mathcal E)$ which is right inverse to the bundle projection.)
\end{Par}

\begin{Prop}
	Let $\mathcal X$ be a principal $\mathcal H$-bundle, and let the super-vector space $V$ carry a linear $\mathcal H$-action. Then $\mathcal X\times^{\mathcal H}\Ber0V\cong\Ber0{\mathcal X\times^{\mathcal H}V}$ as vector bundles. 
\end{Prop}

\begin{Proof}
	Recall the definition of the $\mathcal{GL}(V)$-structure $\mathcal A_V$ on the vector bundle with fibre $V$ associated to the principal bundle $\pi:\mathcal X\to\mathcal B$: $\mathcal A_V(U)$, for any open $U\subset B$, consists of all the $\tilde\tau\in\tau_{\mathcal X\times^{\mathcal H}V}(U)$ induced by the morphisms $\tau':\mathcal X|_{\mathcal U}\times V\to\mathcal U\times V$ defined in \eqref{eq:assochstructure}, where $\tau$ runs through $\tau_{\mathcal X}(U)$. 
	
	Fix an open cover $(U_i)$ of $B$ \scth $\tau_{\mathcal X}(U_i)\neq\vvoid$ \fa $i$, and fix local trivialisations $\tau_i\in\tau_{\mathcal X}(U_i)$. \Fa $i,j$, there exist unique $\vphi_{ij}\in\mathcal{GL}(V)(\mathcal U_{ij})$ \scth $\tau_i'|_{\mathcal X|_{\mathcal U_{ij}}}=\vphi_{ij}.\tau_j'|_{\mathcal X|_{\mathcal U_{ij}}}$, by the proof of \thmref{Prop}{assocbundle}. For the sake of simplicity, we will write this as $\tau_i'=\vphi_{ij}.\tau_j'$. 
	
	If $\tau_i'':\mathcal X|_{\mathcal U_i}\times\Ber0V\to\mathcal U\times\Ber0V$ denotes the morphism defined via \eqref{eq:assochstructure} with $\Ber0V$ as the fibre, then $\tau_i''=\vphi_{ij}.\tau_j''=\Ber0{\vphi_{ij}}.\tau_j''$, by the same proof. Thus, if $\tilde{\tilde\tau}_i\in\mathcal A_{\Ber0V}$ is induced by $\tau_i''$, then $\tilde{\tilde\tau}_i=\Ber0{\vphi_{ij}}.\tilde{\tilde\tau}_j$. This is manifestly the cocycle defining $\Ber0{\mathcal X\times^{\mathcal H}V}$. 
\end{Proof}

\begin{Par}[beraction]
	Let $\mathcal E=(\mathcal X,\mathcal B,p,V)$ be a $\mathcal G$-equivariant vector bundle. Denote by $\mathcal A$ the $\mathcal{GL}(V)$-structure, and by $\tilde\alpha$ and $\alpha$ the actions on $\mathcal X$ and $\mathcal B$, respectively. We use these data to turn $\Ber[_{\mathcal B}]0{\mathcal X}$ into a $\mathcal G$-equivariant vector bundle. This is a somewhat technical construction, but we will need it. 
 
	There exist an open cover $\mathcal U_i$ of $\mathcal B$, and for each $i$, open subsupermanifolds $\mathcal V_k^i$ of $\mathcal G$ and $\mathcal W_k^i$ of $\mathcal B$ together with bundle charts $\tau_i\in\mathcal A(U_i)$ and $\sigma_k^i\in\mathcal A(W_k^i)$, and elements $\gamma_k^i\in\mathcal{GL}(V)(\mathcal V_k^i\times\mathcal W_k^i)$ \scth $\mathcal V_k^i\times\mathcal W_k^i$ form an open cover of $\alpha^{-1}(\mathcal U_i)$, and the following diagrams commute:
	\[
		\xymatrix@C+6ex@R-3ex{%
			\mathcal V_k^i\times\mathcal X|_{\mathcal W_k^i}\ar[r]^-{\id\times\sigma_k^i}\ar[dd]_-{\tilde\alpha}
			&\mathcal V_k^i\times\mathcal W_k^i\times V\ar[rd]^-{(\id,\gamma_k^i)\times\id}\\
			&&\mathcal V_k^i\times\mathcal W_k^i\times\mathcal{GL}(V)\times V\ar[ld]^-{\alpha\times\beta}\\
			\mathcal X|_{\mathcal U_i}\ar[r]_-{\tau_i}&\mathcal U_i\times V
		}
	\]
	Here, we write $\beta$ for the canonical action of $\mathcal{GL}(V)$ on $V$. 
	
	Let $\mathcal U_{ij}=\mathcal U_i\cap\mathcal U_j$, $\mathcal V_{k\ell}^{ij}=\mathcal V_k^i\cap\mathcal V_\ell^j$ and $\mathcal W_{k\ell}^{ij}=\mathcal W_k^i\cap\mathcal W_\ell^j$. There exist elements $\vphi_{ij}\in\mathcal{GL}(\mathcal U_{ij})$ and $\psi_{k\ell}^{ij}\in\mathcal{GL}(\mathcal W_{k\ell}^{ij})$ \scth $\vphi_{ij}.\tau_j=\tau_i$ on $\mathcal X|_{\mathcal U_{ij}}$ and $\psi_{k\ell}^{ij}.\sigma_\ell^j=\sigma_k^i$ on $\mathcal X|_{\mathcal W_{k\ell}^{ij}}$. Since the action $\beta$ is faithful, this implies that $(\vphi_{ij}\circ\alpha)\cdot\gamma_\ell^j=\gamma_k^i\cdot(\psi_{k\ell}^{ij}\circ\pr_2)$ on $\mathcal V_{k\ell}^{ij}\times\mathcal W_{k\ell}^{ij}$ (where the product $\cdot$ of $\mathcal{GL}(V)$-valued morphisms is defined in the obvious way). 
	
	Let $\tilde{\mathcal A}$ denote the $\mathcal{GL}(\Ber0V)$-structure of $\Ber[_{\mathcal B}]0{\mathcal X}$. By its construction, there exist bundle charts $\tilde\tau_i\in\tilde{\mathcal A}(U_i)$ and $\tilde\sigma_k^i\in\tilde{\mathcal A}(W_k^i)$ \scth we have $\Ber0{\vphi_{ij}}.\tilde\tau_j=\tilde\tau_j$ and $\Ber0{\psi_{k\ell}^{ij}}.\tilde\sigma_\ell^j=\tilde\sigma_k^i$. Since 
	\[
		(\Ber0{\vphi_{ij}}\circ\alpha)\cdot\Ber0{\gamma_\ell^j}=\Ber0{\gamma_k^i}\cdot(\Ber0{\psi_{k\ell}^{ij}}\circ\pr_2)\mathtxt{on}\mathcal V_{k\ell}^{ij}\times\mathcal W_{k\ell}^{ij}\ ,
	\]
	there exists a unique morphism $\Ber0{\tilde\alpha}:\mathcal G\times\Ber[_{\mathcal B}]0{\mathcal X}\to\Ber[_{\mathcal B}]0{\mathcal X}$ \scth 
	\[
		\tilde\tau_i\circ\Ber0{\tilde\alpha}=(\alpha\times\tilde\beta)\circ\Parens1{(\id,\Ber0{\gamma_k^i})\times\id}\circ(\id\times\tilde\sigma_k^i)\mathtxt{on}\mathcal V_k^i\times\Ber[_{\mathcal B}]0{\mathcal X}|_{\mathcal W_k^i}
	\]
 	where $\tilde\beta$ denotes the canonical action of $\mathcal{GL}(\Ber0V)$ on $\Ber0V$. It is easy if somewhat tedious to check that $\Ber0{\tilde\alpha}$ is an action, and hence, $\Ber[_{\mathcal B}]0{\mathcal X}$ really is a $\mathcal G$-equivariant vector bundle. 
\end{Par}

\begin{Cor}[quotber]
	There exists an isomorphism $\Ber0{\mathcal G/\mathcal H}\cong\mathcal G\times^{\mathcal H}\Ber0{(\ger g/\ger h)^*}$ of $\mathcal G$-equivariant vector bundles. 
\end{Cor}

\begin{Proof}
	The point is to check the $\mathcal G$-equivariance of the vector bundle isomorphism. By the above considerations, to that end, we need to see that if $\gamma_k^i\in\mathcal{GL}((\ger g/\ger h)^*)(\mathcal V_k^i\times\mathcal W_k^i)$ define the $\mathcal G$-action on $\mathcal G\times^{\mathcal H}(\ger g/\ger h)^*$, then $\Ber0{\gamma_k^i}$ define the $\mathcal G$-action on $\mathcal G\times^{\mathcal H}\Ber0{(\ger g/\ger h)^*}$. By the construction of the vector bundle structure on the associated bundles from the proof of \thmref{Prop}{assocbundle}, the local expressions of the actions factor through a Lie supergroup morphism $\mathcal H\to\mathcal{GL}(V)$ (the same one in both cases), so that the claim follows immediately. 
\end{Proof}

\begin{Par*}
	We come to our first main result.
\end{Par*}

\begin{Th}[unimodequiv]
	Let $\mathcal G$ be a Lie supergroup and $\mathcal H$ a closed subsupergroup with Lie superalgebras $\ger g$ and $\ger h$, respectively. The following are equivalent:
	\begin{enumerate}
		\item The Berezinian bundle $\Ber0{\mathcal G/\mathcal H}$ is $\mathcal G$-equivariantly trivial.
		\item $\Ber0{\mathcal G/\mathcal H}$ has a non-zero $\mathcal G$-invariant global section.
		\item For the action induced by $\Ad^*_\mathcal G$, $(\ger g/\ger h)^*$ is a trivial $\mathcal H$-module. 
	\end{enumerate}
	Here, a global section is a morphism $s:\mathcal G/\mathcal H\to\Pi^q\Ber0{\mathcal G/\mathcal H}$ which is right inverse to the bundle projection. It is called invariant if it is equivariant as a morphism, \emph{i.e.}
	\[
		\Ber0{\alpha_{\mathcal G/\mathcal H}}\circ(\id\times s)=s\circ\alpha\ .
	\]
	Whenever the equivalent conditions are satisfied, then the non-zero $\mathcal G$-invariant section of $\Ber0{\mathcal G/\mathcal H}$ is unique up to constant multiples. 
\end{Th}

\begin{Proof}
	From what we have proved, 1 and 3 are clearly equivalent. The equivalence of 1 and 2 follows by standard procedures. 
\end{Proof}

\begin{Def}
	If the equivalent conditions of \thmref{Th}{unimodequiv} are fulfilled, then the homogeneous supermanifold $\mathcal G/\mathcal H$ is called \emph{unimodular}. (The notion of course depends on the choice of the Lie supergroup $\mathcal G$.)
\end{Def}

\begin{Rem}
	Any Lie supergroup $\mathcal G$ is unimodular as left $\mathcal G$-space. Considered as the quotient $(\mathcal G\times\mathcal G)/\mathcal G$ (for the diagonal action of $\mathcal G$ on $\mathcal G\times\mathcal G$) of the Lie supergroup $\mathcal G\times\mathcal G$, it need not be unimodular. In fact, if $\mathcal G=G$ is a Lie group, then $(G\times G)/G$ is a unimodular $G\times G$-space if and only if $G$ is unimodular as a topological group. This justifies our terminology. 
	
	If $(\mathcal G,\mathcal H)$ is a symmetric pair, \emph{i.e.}~$\mathcal H$ is a closed subsupergroup fixed by an automorphism $\theta$ of $\mathcal G$ of order two, and $H$ is open in $G^\theta$, and $\ger g$ carries a non-degenerate $\mathcal G$-invariant bilinear form for which $\ger h$ is a non-degenerate subspace, then $\mathcal G/\mathcal H$ is unimodular as a $\mathcal G$-space.
\end{Rem}

\section{Fibre integration of Berezinians}\label{sec:berezin}

\subsection{Integration along the fibre}\ \\

\begin{Par}
	If $\mathcal F$ is a sheaf of vector spaces on a topological space $X$, then for $f\in\mathcal F(U)$, the \emph{support} $\supp f\subset X$ is the set of points $x\in X$ where the germ $f_x\neq0$. This set is closed. We let $\Gamma_c(U,\mathcal F)$ be the set of $f\in\mathcal F(U)$ where $\supp f$ is compact. This defines a presheaf $\Gamma_c(\mathcal F)$ on $X$. We call the elements \emph{compactly supported local sections}.
	
	If $\pi:\mathcal X\to\mathcal B$ is a fibre bundle in the category of supermanifolds, then for open subsets $U\subset B$ (!), $\Gamma_{cf}(U,\mathcal O_{\mathcal X})$ denotes the set of $h\in\mathcal O_{\mathcal X}(\pi^{-1}(U))$ \scth $\pi|\supp h:\supp h\to X$ is proper. This defines a presheaf $\Gamma_{cf}(\mathcal O_{\mathcal X})$ on $B$, and its elements are called \emph{compactly supported in the fibre}. One does \emph{not} obtain a presheaf on $X$. Similarly, one may define for any vector bundle on $\mathcal X$ the local sections with compact support or compact support in the fibre (w.r.t.~$\pi$). 
\end{Par}

\begin{Def}
	A supermanifold $\mathcal X$ is called \emph{oriented} if the underlying manifold $X$ is. An isomorphism of oriented supermanifolds is called \emph{orientation preserving} if the underlying isomorphism of oriented manifolds is orientation preserving. 
	
	If $\mathcal X$ is an oriented supermanifold, then there exists a linear morphism of presheaves on $X$ \cite[\S~6.2]{manin}
	\[
		\int_{\mathcal X}:\Gamma_c(\Ber0{\mathcal X})\to\reals_X\ ,
	\]
	which in the case of a manifold $\mathcal X=X$ is the integration of volume forms. Moreover, if $\vphi:\mathcal X\to\mathcal Y$ is an orientation preserving isomorphism of oriented supermanifolds, then \cite[Theorem 2.4.5]{leites}
	\begin{equation}\label{eq:trans}
		\int_{\mathcal X}\vphi^*\omega=\int_{\mathcal Y}\omega\mathfa\omega\in\Gamma_c(U,\Ber0{\mathcal Y})
	\end{equation}
	where we recall the definition of $\vphi^*\omega$ from \thmref{Par}{berpbdef}. 
	
	Let $\mathcal E=(\mathcal X,\mathcal B,\pi,\mathcal F)$ be a fibre bundle. It is called \emph{oriented} if $\mathcal X$, $\mathcal B$ and $\mathcal F$ are oriented and it is supplied with an \emph{oriented bundle atlas}. The latter is the data of an open covering $(\mathcal U_i)$ of $\mathcal B$, and of local trivialisations $\tau_i\in\tau_{\mathcal X}(U_i)$ which are \emph{orientation preserving} isomophisms $\tau_i:\mathcal X|_{\mathcal U_i}\to\mathcal U_i\times\mathcal F$ where on the right hand side, we take the product orientation. 
\end{Def}

\begin{Par}
	We wish to define fibre integration of Berezinians. To that end, we have to introduce topologies on the presheaves $\mathcal O_{\mathcal X}$ and $\Gamma_c(\mathcal O_{\mathcal X})$ for any supermanifold $\mathcal X$, and on $\Gamma_{cf}(\mathcal O_{\mathcal X\times\mathcal Y})$ for any direct product of supermanifolds. 
	
	Let $U\subset X$ be open. For any compact $K\subset U$, and any differential operator $D\in\mathcal D_{\mathcal X}(U)$, we define a seminorm $p_{K,D}$ on $\mathcal O_{\mathcal X}(U)$ as follows:
	\[
		p_{K,D}(h)=\sup\nolimits_{x\in K}\Abs0{(Dh)(x)}\ .
	\]
	
	Here, the value $h(x)$ at $x$ of $h\in\mathcal O_{\mathcal X}(U)$ is defined as the value at $x$ of the image of $h$ in $\mathcal C^\infty_X(U)$. 
\end{Par}

\begin{Lem}
	Let $\mathcal X$ be a supermanifold and $U\subset X$ an open subset. Then $\mathcal O_{\mathcal X}(U)$ is a nuclear and $m$-convex Fr\'echet algebra. (Where we recall that a topological algebra is called \emph{$m$-convex} if its topology is the locally convex topology defined by a family of submultiplicative seminorms.)
\end{Lem}

\begin{Proof}
	If $\mathcal U$, the open subsupermanifold of $\mathcal X$ with base $U\!$, is a superdomain, then $\mathcal O_{\mathcal X}(U)\cong\Ct[^\infty]0{U}\otimes\bigwedge(\reals^q)$ as topological algebras (where $\dim\mathcal U=p|q$). The right hand side is certainly an $m$-convex Fr\'echet algebra, and both tensor factors are nuclear. Thus, in this case, the statement follows from the fact that the projective tensor product of nuclear spaces is nuclear \cite[Chapter II, \S 7.5]{schaefer-tvs}. 

	In general, $\mathcal U$ has an open cover $(\mathcal U_i)$ by superdomains. Then $\mathcal O_{\mathcal X}(U)$ carries the initial locally convex topology with respect to the restriction maps $\vrho^U_{U_i}:\mathcal O_{\mathcal X}(U)\to\mathcal O_{\mathcal X}(U_i)$. Thus, $\mathcal O_{\mathcal X}(U)$ is a projective limit of complete nuclear locally convex spaces and as such, complete and nuclear \cite[Chapter II, \S~5.3; Chapter III, \S~7.4, Corollary]{schaefer-tvs}. It follows also that it is $m$-convex. Since $X$ is second countable, the topology of $\mathcal O_{\mathcal X}(U)$ is generated by a countable family of seminorms. Therefore, $\mathcal O_{\mathcal X}(U)$ is metrisable, and hence, a Fr\'echet space.
\end{Proof}

\begin{Lem}[tensortop]
	Let $\mathcal X\times\mathcal Y$ be supermanifolds, $U\subset X$, $V\subset V$ be open subsets. If $\tau$ is any tensor product topology, then $\mathcal O_{\mathcal X}(U)\,\Hat\otimes_\tau\,\mathcal O_{\mathcal Y}(V)\cong\mathcal O_{\mathcal X\times\mathcal Y}(U\times V)$ via the continuous linear extension of the map $\phi:g\otimes h\mapsto\pr_1^*g\cdot\pr_2^*h$. 
\end{Lem}

\begin{Proof}
	The map $\phi:\mathcal O_{\mathcal X}(U)\otimes\mathcal O_{\mathcal Y}(V)\to\mathcal O_{\mathcal X\times\mathcal Y}(U\times V)$ is injective and has dense image. This is easily checked for superdomains, and the general statement follows by a projective limit argument, or by considering partitions of unity.

	Both tensor factors are nuclear Fr\'echet spaces, so the point is to show that the restriction $\psi$ of $\phi$ to $\mathcal O_{\mathcal X}(U)\times\mathcal O_{\mathcal Y}(V)$ is a continuous bilinear map. We need to show that $p_{K,D}\circ\psi$ is continuous for any compact $K\subset U\times V$ and any $D\in\mathcal D_{\mathcal X\times\mathcal Y}(U\times V)$. Since $K$ is compact, it has a finite cover by superdomains, and we may restrict our attention to the case where $\mathcal U$ and $\mathcal V$ are superdomains. By an easy estimate, the case of a general differential operator $D$ is reduced to the case of one which is a polynomial in the vector fields $\frac\partial{\partial x_i}$, $\frac\partial{\partial\xi_j}$, $\frac\partial{\partial y_k}$, $\frac\partial{\partial\eta_\ell}$ where $(x_i,\xi^j)$, $(y_k,\eta^\ell)$ are coordinates for $\mathcal U$ and $\mathcal V$, respectively. But this case is entirely trivial. 
\end{Proof}

\begin{Par}[tensortop]
	In what follows, we denote by $\Hat\otimes$ the completed (and graded) \emph{projective} tensor product of locally convex vector spaces.  
	
	Let $\Gamma_c(U,\mathcal O_{\mathcal X})$ be topologised as the inductive limit with respect to the inclusions $\Gamma_K(U,\mathcal O_{\mathcal X})\subset\Gamma_c(U,\mathcal O_{\mathcal X})$, where for any compact $K\subset U$, $\Gamma_K(U,\mathcal O_{\mathcal X})$ is the subspace of $\mathcal O_{\mathcal X}(U)$ consisting of all $h$ with $\supp h\subset K$. Since $U$ is second countable, the topology is defined by a \emph{countable} inductive limit. Thus, $\Gamma_c(U,\mathcal O_{\mathcal X})$ is an $(LF)$-space, and it is nuclear \cite[Chapter III, \S~7.4, Corollary]{schaefer-tvs}.
	
	Next, we fix a second supermanifold $\mathcal Y$. Consider the first projection $p:\mathcal X\times\mathcal Y\to\mathcal X$ and $\Gamma_{cf}(U,\mathcal O_{\mathcal X\times\mathcal Y})$. This space is topologised as the inductive limit with respect to the inclusions $\Gamma_{cf,K}(U,\mathcal O_{\mathcal X\times\mathcal Y})\subset\Gamma_{cf}(U,\mathcal O_{\mathcal X\times\mathcal Y})$ where for any closed $K\subset X\times Y$ with $p|_K:K\to X$ proper, the former space is defined to be the subspace of $\mathcal O_{\mathcal X\times\mathcal Y}$ which consists of those $h$ \scth $\supp h\subset K$. This makes $\Gamma_{cf}(U,\mathcal O_{\mathcal X\times\mathcal Y})$ a nuclear $(LF)$-space. Moreover, we have as above that $\Gamma_{cf}(U,\mathcal O_{\mathcal X\times\mathcal Y})=\mathcal O_{\mathcal X}(U)\,\Hat\otimes\,\Gamma_c(Y,\mathcal O_{\mathcal Y})$. 
	
	If we are given a vector bundle $\mathcal E$ on $\mathcal X$, then the sheaf of sections $\Gamma(\mathcal E)$ is locally free over $\mathcal O_{\mathcal X}$. On trivialising open subsupermanifolds of $\mathcal X$, this determines a topology on $\Gamma(\mathcal E)$ (the product topology). If $U\subset X$ is any open subset, then we let $\Gamma(U,\mathcal E)$ be equipped with the initial locally convex topology with respect to all restrictions to trivialising open subsets (by the open mapping theorem, this is well-defined). 
	
	Finally, one defines topologies on $\Gamma_c(\mathcal E)$ and $\Gamma_{cf}(\mathcal E)$ (in the case of a vector bundle over a direct product of supermanifolds) in the same way as above, namely, by taking inductive limits over subsets $K$ which are compact or closed and \scth $p|_K:K\to X$ is proper, respectively. The statements for projective tensor products carry over. 
\end{Par}

\begin{Prop}[fibreint]
	Let $\mathcal E=(\mathcal X,\mathcal B,\pi,\mathcal F)$ be an oriented fibre bundle in supermanifolds where $\dim\mathcal B=m|n$ and $\dim\mathcal F=p|q$. There is an even morphism $\pi_!:\Gamma_{cf}(\Ber0{\mathcal X})\to\Gamma(\Ber0{\mathcal B})$ of graded presheaves over $B$, \scth
	\begin{equation}\label{eq:globfibintsupp}
		\pi_!(\pi^*h\cdot\omega)=h\cdot\pi_!(\omega)\nd\supp\pi_!(\omega)\subset\pi(\supp\omega)
	\end{equation}
	\fa open $U\subset B$, $h\in\mathcal O_{\mathcal B}(U)$, $\omega\in\Gamma_{cf}(U,\Ber0{\mathcal X})$, and 
	\begin{equation}\label{eq:fibreintid}
		\int_{\mathcal X}\omega=(-1)^{(m+n)q}\cdot\int_{\mathcal B}\pi_!(\omega)
	\end{equation}
	\fa $\omega\in\Gamma_c(\pi^{-1}(U),\Ber0{\mathcal X})$. 
	
	Moreover, if $\mathcal E'$ is another fibre bundle over $\mathcal B$ and $\vphi:\mathcal X'\to\mathcal X$ is an oriented isomorphism of the total spaces \scth $\pi'=\pi\circ\vphi$, then 
	\begin{equation}\label{eq:fibintiso}
		\pi'_!(\vphi^*\omega)=\pi_!(\omega)\mathfa\omega\in\Gamma_{cf}(U,\Ber0{\mathcal X})\ .
	\end{equation}
\end{Prop}

In the \emph{proof}, we first establish the local picture. 

\begin{Lem}[fibreint]
	Let $\mathcal X$ and $\mathcal Y$ be supermanifolds where $\mathcal Y$ is oriented, and set $p=\pr_1:\mathcal X\times\mathcal Y\to\mathcal X$. Let $\mathcal E$ be a vector bundle on $\mathcal X$. For any open $U\subset X$, we define the map $p_!:\Gamma(U,\mathcal E)\otimes\Gamma_c(Y,\Ber0{\mathcal Y})\to\Gamma(U,\mathcal E)$ by 
	\[
	p_!(\omega_1\otimes\omega_2)=\omega_1\cdot\int_{\mathcal Y}\omega_2\ .
	\]
	Then $p_!$ has an extension $\Gamma_{cf}(U,\mathcal E\boxtimes\Ber0{\mathcal Y})\to\Gamma(U,\mathcal E)$ uniquely determined by the requirement that it be continuous and linear (here, $\boxtimes$ denotes the graded external tensor product). Thus, we obtain an even morphism  of graded presheaves on $X$, $p_!:\Gamma_{cf}(\mathcal E\boxtimes\Ber0{\mathcal Y})\to\Gamma(\mathcal E)$, which satisfies
	\begin{equation}\label{eq:locfibintsupp}
	p_!(p^*h\cdot\omega)=h\cdot p_!(\omega)\nd
	\supp p_!(\omega)\subset p(\supp\omega)
	\end{equation}
	\fa $h\in\mathcal O_{\mathcal X}(U)$, $\omega\in\Gamma_{cf}(U,\mathcal E\boxtimes\Ber0{\mathcal Y})$.
\end{Lem}

\begin{Proof}
	From our remarks in \thmref{Par}{tensortop} on projective tensor products, it follows that $\Gamma_{cf}(U,\mathcal E\boxtimes\Ber0{\mathcal Y})=\Gamma(U,\mathcal E)\,\Hat\otimes\,\Gamma_c(Y,\Ber0{\mathcal Y})$. Thus, to prove unique existence of a continuous linear extension, it suffices to see that for any open subset $U\subset X$, $p_!$ restricts to a continuous linear bilinear map on $\Gamma(U,\mathcal E)\times\Gamma_c(Y,\Ber0{\mathcal Y})$. Since $\int_{\mathcal Y}$ is continuous on any $\Gamma_K(Y,\Ber0{\mathcal Y})$ for $K\subset Y$ compact, the assertion follows immediately. 
	
	To check \eqref{eq:locfibintsupp}, we note that the second statement follows from the first; and to prove the first, it suffices to consider $p_!$ on the algebraic tensor product. Here, we have
	\begin{multline*}
		p_!\Parens1{p^*h\cdot(\omega_1\otimes\omega_2)}=p_!\Parens1{(h\cdot\omega_1)\otimes\omega_2}
		=h\cdot\omega_1\cdot\int_{\mathcal Y}\omega_2=h\cdot p_!(\omega_1\otimes\omega_2)\ .
	\end{multline*}
	This proves the assertion.
\end{Proof}

\noindent\textbf{Proof of \thmref{Prop}{fibreint}.$\quad$}First, note that $\Ber0{\mathcal Y\times\mathcal Z}=\Ber0{\mathcal Y}\boxtimes\Ber0{\mathcal Z}$. Indeed,  for super-vector spaces $V$ and $W$, it follows from the definition of $\Ber0V$ that there is a canonical isomorphism $\Ber0{V\oplus W}=\Ber0V\otimes\Ber0W$. Thus, \thmref{Lem}{fibreint} gives a morphism $p_!:\Gamma_{cf}(\Ber0{\mathcal Y\times\mathcal Z})\to\Gamma(\Ber0{\mathcal Y})$. 

	Let $(\mathcal U_i)$, $\tau_i\in\tau_{\mathcal X}(U_i)$ be the data of an oriented bundle atlas where we assume that $U_i$ is relatively compact in $B$ for any $i$, and set $\sigma_i=\tau_i^{-1}$. Let $\vphi_i\in\Ct[^\infty]0B$ form a partition of unity on $B$ subordinate to $(U_i)$. We write $p$ for the first projections $\mathcal U_i\times\mathcal F\to\mathcal U_i$. 
	
	For any open $U\subset B$, we define
	\[
		\pi_!(\omega)=\textstyle\sum_ip_!(p^*\vphi_i\cdot\sigma_i^*(\omega))\mathfa\omega\in\Gamma_{cf}(U,\Ber0{\mathcal X})\ .
	\]

	As usual, the issue is to see that this is well-defined independently of all choices. First, to prove independence of the partition of unity, one applies \eqref{eq:locfibintsupp}. Next, the independence of the definition on the choice of oriented bundle atlas follows by considering elementary tensors and applying \eqref{eq:trans}. Hence, the definition is independent of all choices, and in particular, it follows that $\pi_!$ is linear and defines a morphism of presheaves. Then \eqref{eq:globfibintsupp} follows from \eqref{eq:locfibintsupp}, and \eqref{eq:fibintiso} is also an application of \eqref{eq:trans}.  
	
	Finally, to prove \eqref{eq:fibreintid}, we may assume that $\mathcal X=\mathcal B\times\mathcal F$ and that $\pi$ is the first projection, so we are in the situation of \thmref{Lem}{fibreint}. By shrinking coordinate charts in the first factor, and by taking partitions of unity in the second, we may assume that the line bundles $\Ber0{\mathcal B}$ and $\Ber0{\mathcal F}$ are trivial, and that $\mathcal B$ and $\mathcal F$ are superdomains. By \eqref{eq:globfibintsupp}, it suffices to prove 
	\[
		\int_{\mathcal B\times\mathcal F}(\omega_1\otimes\omega_2)=(-1)^{(m+n)q}\cdot\int_{\mathcal B}\omega_1\cdot\int_{\mathcal F}\omega_2\tag{$*$}
	\]
	\fa $\omega_1\in\Gamma_c(B,\Ber0{\mathcal B})$ and $\omega_2\in\Gamma_c(F,\Ber0{\mathcal F})$. 
	
	Fix coordinates $(x_1,\dotsc,x_m,\xi^1,\dotsc,\xi^n)$ of $\mathcal B$ and $(y_1,\dotsc,y_p,\eta^1,\dotsc,\eta^q)$ of $\mathcal F$. We have global sections $D(x,\xi)$ and $D(y,\eta)$ of $\Ber0{\mathcal B}$ and $\Ber0{\mathcal F}$, respectively \cite[Chapter 3, \S~4.7]{manin}. Here, we recall that 
	\[
		\Pi^m(D(x,\xi))\equiv\Pi(dx_1)\dotsm\Pi(dx_m)\cdot\tfrac\partial{\partial\xi^1}\dotsm\tfrac\partial{\partial\xi^n}
	\]
	(modulo the image of $\Pi$). 
	Moreover, $D(x,y,\xi,\eta)=(-1)^{np}D(x,\xi)\otimes D(y,\eta)$ is a global section of $\Ber0{\mathcal B\times\mathcal F}$. 
	
	Let $g\in\mathcal O_{\mathcal B}(B)$, $h\in\mathcal O_{\mathcal F}(F)$. Write $g=\sum_\alpha\xi^\alpha g_\alpha$ and $h=\sum_\beta\eta^\beta h_\beta$ where $g_\alpha$, $h_\beta$ are even, and $dx=dx_1\dotsm dx_m$ and $dy=dy_1\dotsc dy_q$. The Berezin integral on $\mathcal B$ is given by $\int_{\mathcal B}D(x,\xi)\cdot g=(-1)^{mn}\int_Bg_{1,\dotsc,1}\,dx$. Thus 
	\begin{align*}
		\int_{\mathcal B\times\mathcal F}D(x,y,\xi,\eta)\,(g\otimes h)&=(-1)^{(m+p)(n+q)}\cdot\int_{B\times F}g_{1,\dotsc,1}h_{1,\dotsc,1}\,dxdy\\
		&=(-1)^{(m+p)(n+q)}\cdot\int_Bg_{1,\dotsc,1}\,dx\cdot\int_Fh_{1,\dotsc,1}\,dy\\
		&=(-1)^{mq+np}\cdot\int_{\mathcal B}D(x,\xi)\cdot g\cdot\int_{\mathcal F}D(y,\eta)\cdot h
	\end{align*}
	Finally, note 
	\[
	(D(x,\xi)\cdot g)\otimes(D(y,\eta)\cdot h)=(-1)^{n\Abs0g}(D(x,\xi)\otimes D(y,\eta))\cdot(g\otimes h)
	\]
	and $\Abs0{\xi^{1,\dotsc,1}g_{1,\dotsc,1}}=q$. This proves ($*$), and therefore, the assertion. 
\qed\ \\ 

\subsection{`Fubini' formula for quotients $\mathcal G/\mathcal H$}\ \\

\begin{Par}
	Let $\mathcal G$ be a Lie supergroup and $\mathcal H$ be a closed Lie subsupergroup of $\mathcal G$ \scth $\mathcal G/\mathcal H$ is a unimodular $\mathcal G$-space. By \thmref{Th}{unimodequiv}, there exist non-zero $\mathcal G$-invariant section $\omega_{\mathcal G}$ of $\Ber0{\mathcal G}$ and $\omega_{\mathcal G/\mathcal H}$ of $\Ber0{\mathcal G/\mathcal H}$, and a non-zero $\mathcal H$-invariant section $\omega_{\mathcal H}$ of $\Ber0{\mathcal H}$. For $f\in\Gamma_c(G,\mathcal O_{\mathcal G})$, we define
	\[
		\int_{\mathcal G} f=\int_{\mathcal G}f\cdot\omega_{\mathcal G}\ ,
	\]
	and similarly for $\mathcal H$ and $\mathcal G/\mathcal H$. 
\end{Par}

\begin{Prop}[fubformula]
	Retain the above assumptions. For a suitable normalisation of $\omega_{\mathcal G}$, $\omega_{\mathcal H}$ and $\omega_{\mathcal G/\mathcal H}$, we have, for each bundle chart $\tau:\mathcal G|_{\mathcal U}\to\mathcal U\times\mathcal H$ of the principle $\mathcal H$-bundle $\mathcal G\to\mathcal G/\mathcal H$, $\omega_{\mathcal G}=\tau^*(\omega_{\mathcal G/\mathcal H}\otimes\omega_{\mathcal H})$. 
\end{Prop}

\begin{Proof}
	Let $\mathcal A$ denote the $\mathcal H$-structure of $\pi:\mathcal G\to\mathcal G/\mathcal H$. For any open subset $U\subset G/H$ \scth $\mathcal A(U)\neq\vvoid$ and any $\tau\in\mathcal A(U)$, we may define a local section $\omega_\tau\in\Gamma(\pi^{-1}(U),\Ber0{\mathcal G})$ by $\omega_\tau=\tau^*(\omega_{\mathcal G/\mathcal H}\otimes\omega_{\mathcal H})$. Also, certainly, $\omega_\tau$ is non-zero on $\mathcal G|_{\mathcal U}$. 
	
	Let $V\subset G/H$ be an open subset \scth $\mathcal A(V)\neq\vvoid$ and $U\cap V\neq\vvoid$. Let $\sigma\in\mathcal A(V)$. There exists a unique $\vphi\in\mathcal H(\mathcal U\cap\mathcal V)$ \scth $\tau=\vphi.\sigma$. Then 	
	\[
		(\omega_{\mathcal G/\mathcal H}\circ\pr_1)\otimes\Parens1{\Ber0{m_{\mathcal H}}\circ(\id\times\omega_{\mathcal H})\circ(\pr_2,\pr_3)}=(\omega_{\mathcal G/\mathcal H}\otimes\omega_{\mathcal H})\circ(\id\times m_{\mathcal H})
	\]
	by the invariance of $\omega_{\mathcal H}$. It follows that on $\mathcal U\cap\mathcal V$, 	
	\[
		(\omega_{\mathcal G/\mathcal H}\otimes\omega_{\mathcal H})\circ(\pr_1,m_{\mathcal H}\circ(\vphi\times\id))=\Parens1{(\omega_{\mathcal G/\mathcal H}\circ\pr_1)\otimes\Ber0{m_{\mathcal H}}}\circ(\pr_1,\vphi\circ\pr_1,\omega_{\mathcal H}\circ\pr_2)\ ,
	\]
	so 
	\[
		\omega_\tau=\sigma^*(\pr_1,m_{\mathcal H}\circ(\vphi\times\id))^*(\omega_{\mathcal G/\mathcal H}\otimes\omega_{\mathcal H})=\omega_\sigma\ .
	\]
	Hence, there exists a unique $\omega\in\Gamma(G,\Ber0{\mathcal G})$ \scth $\omega|_{\pi^{-1}(U)}=\omega_\tau$ (where $\tau\in\mathcal A(U)$ is arbitrary), for all open sets $U\subset G/H$. 
	
	To prove the assertion, we have to see that $\omega$ is $\mathcal G$-invariant. Choose open subsets $U\subset G$, $V,W\subset G/H$ \scth $m_{\mathcal G}(U\times\pi^{-1}(V))\subset\pi^{-1}(W)$ and $\mathcal A(V),\mathcal A(W)\neq\vvoid$. Let $\tau\in\mathcal A(V)$ and $\sigma\in\mathcal A(W)$. By the proof of \thmref{Prop}{principalbundle}, there exist sections $t:\mathcal V\to\mathcal G|_{\mathcal V}$ and $s:\mathcal W\to\mathcal G|_{\mathcal W}$ of $\pi$ \scth $\tau^{-1}=m_{\mathcal G}\circ(t\times\id)$ and $\sigma^{-1}=m_{\mathcal G}\circ(s\times\id)$. 
	
	Define $\vphi:\mathcal U\times\mathcal V\to\mathcal G$ by $\vphi=m_{\mathcal G}\circ(i_{\mathcal G}\circ s\circ\alpha,m_{\mathcal G}\circ(\id\times t))$. Here, $\alpha$ denotes the action of $\mathcal G$ on $\mathcal G/\mathcal H$, and $i_{\mathcal G}$ denotes the inversion of $\mathcal G$.
	
	We have $\pi\circ m_{\mathcal G}=\alpha\circ(\id\times\pi)$, so 
	\[
		\pi\circ s\circ\alpha=\alpha=\pi\circ m_{\mathcal G}\circ(\id\times t)\ ,
	\]
	so that $s\circ\alpha$ and $m_{\mathcal G}\circ(\id\times t)$ induce a morphism $\tilde\vphi:\mathcal U\times\mathcal V\to\mathcal G\times_{\mathcal G/\mathcal H}\mathcal G$. By the effectiveness of the quotient $\mathcal G/\mathcal H$ (\thmref{Th}{quotexist}), $\iota=(m_{\mathcal G},\pr_1):\mathcal G\times\mathcal H\to\mathcal G\times\mathcal G$ induces by corestriction an isomorphism $\mathcal G\times\mathcal H\to\mathcal G\times_{\mathcal G/\mathcal H}\mathcal G$. Hence, there exists a morphism $\psi:\mathcal U\times\mathcal V\to\mathcal H$ \scth $s\circ\alpha=m_{\mathcal G}\circ(m_{\mathcal G}\circ(\id\times t),\psi)$ (take $\psi=\pr_2\circ\,\iota^{-1}\circ\tilde\vphi$). An easy computation shows that $\vphi=i_{\mathcal G}\circ\psi$, so that the corestriction $\vphi:\mathcal U\times\mathcal V\to\mathcal H$ exists. Moreover, the following diagram commutes:
	\[
		\xymatrix@C+6ex{%
			\mathcal U\times\mathcal G|_{\mathcal V}\ar[rr]^-{m_{\mathcal G}}
			&&\mathcal G|_{\mathcal W}\\
			\mathcal U\times\mathcal V\times\mathcal H\ar[u]^-{\id\times\tau^{-1}}
			\ar[rd]_-{(\id,\vphi)\times\id\qquad}
			&&\mathcal W\times\mathcal H\ar[u]_-{\sigma^{-1}}\\
			&\mathcal U\times\mathcal V\times\mathcal H\times\mathcal H\ar[ru]_-{\quad\alpha\times m_{\mathcal H}}
		}
	\]
	Indeed, if $m_{\mathcal G}^{(2)}=m_{\mathcal G}\circ(m_{\mathcal G}\times\id)=m_{\mathcal G}\circ(\id\times m_{\mathcal G})$, then 
	\begin{align*}
		\sigma^{-1}\circ(\alpha&\times m_{\mathcal H})\circ((\id,\vphi)\times\id)\\
		&=m_{\mathcal G}\circ(s\times\id)\circ(\alpha\times m_{\mathcal H})\circ((\id,\vphi)\times\id)\\
		&=m_{\mathcal G}^{(2)}\circ(s\circ\alpha\circ(\pr_1,\pr_2),i_{\mathcal G}\circ s\circ\alpha\circ(\pr_1,\pr_2),m_{\mathcal G}^{(2)}\circ(\id\times t\times\id))\\
		&=m_{\mathcal G}^{(2)}\circ(\id\times t\times\id)=m_{\mathcal G}\circ(\id\times\tau^{-1})\ .
	\end{align*}
	 
	By the invariance of $\omega_{\mathcal G/\mathcal H}$ and $\omega_{\mathcal H}$, 
	\[
		(\omega_{\mathcal G/\mathcal H}\otimes\omega_{\mathcal H})\circ(\alpha\times m_{\mathcal H})
		=\Ber0{\alpha\times m_{\mathcal H}}\circ(\id\times(\omega_{\mathcal G/\mathcal H}\otimes\omega_{\mathcal H}))
	\]
	on $\mathcal U\times\mathcal V\times\mathcal H$. Hence, a similar computation as for the well-definedness of $\omega$ shows the relation $\Ber0{m_{\mathcal G}}\circ(\id\times\omega)=\omega\circ m_{\mathcal G}$ on $\mathcal U\times\mathcal G|_{\mathcal V}$. This proves the assertion. 
\end{Proof}

\begin{Par}
	For any $f\in\Gamma_c(G,\mathcal O_{\mathcal G})$, there exists a unique $h\in\Gamma_c(G/H,\mathcal O_{\mathcal G/\mathcal H})$ \scth $\pi_!(f\cdot\omega_{\mathcal G})=h\cdot\omega_{\mathcal G/\mathcal H}$. We write $f_{\mathcal H}=h$. Then 
	\begin{equation}\label{eq:fibreint}
		\int_{\mathcal G}f=(-1)^{\dim\ger h_1\cdot\dim\ger g/\ger h}\cdot\int_{\mathcal G/\mathcal H}f_{\mathcal H}\mathfa f\in\Gamma_c(G,\mathcal O_{\mathcal G})\ .
	\end{equation}
\end{Par}

\begin{Cor}[fubformula]
	Retain the assumptions of \thmref{Prop}{fubformula}. For a suitable normalisation of the Berezinians, for all bundle charts $\tau:\mathcal U\times\mathcal H\to\mathcal G|_{\mathcal U}$, and all $f\in\Gamma_c(G,\mathcal O_{\mathcal G})$, the following identity holds:
	\begin{equation}\label{eq:fubformula}
	f_{\mathcal H}|_U=p_!(\tau^*f\cdot(1\otimes\omega_{\mathcal H}))
	\end{equation}
	Here, $p:\mathcal U\times\mathcal H\to\mathcal U$ is the first projection, and we denote by $p_!$ the linear map $\Gamma_{cf}(U,\mathcal O_{\mathcal G/\mathcal H}\otimes\Ber0{\mathcal H})\to\mathcal O_{\mathcal G/\mathcal H}(U)$ defined in \thmref{Lem}{fibreint} for the case of the trivial bundle on $\mathcal G/\mathcal H$. 
\end{Cor}

\begin{Proof}
	This follows from \thmref{Prop}{fibreint}, the definitions in \thmref{Lem}{fibreint}, \thmref{Prop}{fubformula}, and standard considerations using tensor products. 
\end{Proof}

\begin{Rem}
	The formulae \eqref{eq:fibreint} and \eqref{eq:fubformula} constitute the graded version of the classical `Fubini' formula 
	\[
		\int_G f(g)\,dg=\int_{G/H}\Parens3{\int_H f(gh)\,dh}\,d\dot g\mathfa f\in\Cc0G
	\]
	valid for any unimodular homogeneous $G$-space $G/H$. 
\end{Rem}

\subsection{Products of subsupergroups}\ \\

In this section, we will generalise the results \cite[Chapter I, Lemma 1.11, Proposition 1.12]{helgason-gaga} to the setting of homogeneous supermanifolds.\ \\

\begin{Par}[berchart]
	Let $V$ be a finite-dimensional super-vector space, $x_1,\dotsc,x_n$ a homogeneous basis, and $\xi_1,\dotsc,\xi_n$ its dual basis. We extend the definition of $D(x,\xi)$ for graded bases to this setting. Let $D(x)=D(x_1,\dotsc,x_n)\in\Ber0V$ be the element which is represented modulo the image of the differential $\Pi$ on $S(\Pi V\oplus V^*)$ by 
	\[
		y_1\dotsc y_n\mathtxt{where}y_i=\begin{cases}\Pi(x_i)&\Abs0{x_i}=0\\\xi_i&\Abs0{x_i}=1\end{cases}
	\]
	
	We define a \emph{canonical} perfect pairing $\Dual0\cdot\cdot:\Ber0{V^*}\otimes\Ber0{V}\to\reals$ by
	\[
	\Dual0{D(\xi_n,\dotsc,\xi_1)}{D(x_1,\dotsc,x_n)}=1\ .
	\]
	The definition is independent of the choice of bases \cite[Lemma 1.4]{duflo-petracci}.
	
	Next, recall from \thmref{Prop}{quottangent} that for any Lie supergroup $\mathcal G$ and any closed subsupergroup $\mathcal H$, we have a $\mathcal G$-equivariant vector bundle isomorphism $T(\mathcal G/\mathcal H)\cong\mathcal G\times^{\mathcal H}\ger g/\ger h$. In particular, for any bundle chart $\tau:\mathcal U\times\mathcal H\to\mathcal G|_{\mathcal U}$, there is a bundle chart $\mathcal U\times\ger g/\ger h\to T(\mathcal G/\mathcal H)|_{\mathcal U}$ which we denote by $T\tau$. 
	
	Explicitly, it is given as follows. If $t:\mathcal U\to\mathcal G$ is a local section of $\pi:\mathcal G\to\mathcal G/\mathcal H$ \scth $\tau=m_{\mathcal G}\circ(t\times\id)$, then $T\tau=\alpha_{\mathcal G/\mathcal H}\circ(t\times\eta)$ where $\alpha_{\mathcal G/\mathcal H}:\mathcal G\times\mathcal G/\mathcal H\to\mathcal G/\mathcal H$ is the $\mathcal G$-action, and $\eta:\ger g/\ger h\to T(\mathcal G/\mathcal H)$ is (again) the canonical inclusion of the fibre at the base point $o=H\in G/H$. 
	
	Since we also have an induced $\mathcal G$-equivariant vector bundle isomorphism $\Ber0{\mathcal G/\mathcal H}\cong\mathcal G\times^{\mathcal H}\Ber0{(\ger g/\ger h)^*}$, by \thmref{Cor}{quotber}, the construction of bundle charts carries over to this situation. Indeed, associated with $\tau$, there is a bundle chart $\Ber0\tau:\mathcal U\times\Ber0{(\ger g/\ger h)^*}\to\Ber0{\mathcal G/\mathcal H}|_{\mathcal U}$, and it is given explicitly by $\Ber0\tau=\Ber0{\alpha_{\mathcal G/\mathcal H}}\circ(t\times\eta)$ where $\eta:\Ber0{(\ger g/\ger h)^*}\to\Ber0{\mathcal G/\mathcal H}$ is the canonical inclusion of the fibre at the base point $o=H\in G/H$, and $\Ber0{\alpha_{\mathcal G/\mathcal H}}$ is the Berezinian action defined in \thmref{Par}{beraction}. 
\end{Par}

\begin{Lem}[berhomsplociso]
	Let $\mathcal G/\mathcal H$ and $\mathcal S/\mathcal T$ be unimodular as a $\mathcal G$- and as an $\mathcal S$-space, respectively. Fix invariant forms $\omega_{\mathcal G/\mathcal H}$, $\omega_{\mathcal S/\mathcal T}$, and base points $o=H\in G/H$ and $o'=T\in S/T$. Assume that $\mathcal G/\mathcal H$ and $\mathcal S/\mathcal T$ have common graded dimension $p|q$, $n=p+q$, that $\phi:\mathcal G/\mathcal H\to\mathcal S/\mathcal T$ is a local isomorphism onto an open subsupermanifold of $\mathcal S/\mathcal T$, and that $\phi(o)=o'$. 
	
	For any bundle charts $\tau:\mathcal U\times\mathcal H\to\mathcal G|_{\mathcal U}$ and $\sigma:\mathcal V\times\mathcal T\to\mathcal S|_{\mathcal V}$ \scth $\phi:\mathcal U\to\mathcal V$ is an open embedding, there exists $d\in\mathcal O_{\mathcal G/\mathcal U}(U)$ \scth 
	\[
	\phi^*(\omega_{\mathcal S/\mathcal T})=d\cdot\omega_{\mathcal G/\mathcal H}\mathtxt{on}\mathcal U\ .
	\]
	
	The superfunction $d$ may be computed as follows: There are local sections $t:\mathcal U\to\mathcal G$ and $s:\mathcal V\to\mathcal H$ \scth $\tau=m_{\mathcal G}\circ(t\times\id)$ and $\sigma=m_{\mathcal S}\circ(s\times\id)$. Let $x_1,\dotsc,x_n$, $y_1,\dotsc,y_n$ be homogeneous bases of $\ger g/\ger h$ and $\ger s/\ger t$, respectively, such that
	\[
		\Dual0{(\omega_{\mathcal G/\mathcal H})_o}{D(x_1,\dotsc,x_n)}=\Dual0{(\omega_{\mathcal S/\mathcal T})_{o'}}{D(y_1,\dotsc,y_n)}\ .
	\]
	 Consider the local expression $\vphi:\mathcal U\to\Hom0{\ger g/\ger h,\ger s/\ger t}_0$ of $T\phi$; according to \eqref{eq:vbmordef}, 
	\begin{equation}\label{eq:locexpr}
	T\phi\circ T\tau=T\sigma\circ(\id\times\eps)\circ((\id,\vphi)\times\id)
	\end{equation}
	where $T\tau$, $T\sigma$ where defined in \thmref{Par}{berchart}. Define $\xi=(\xi_{ij})\in\GL\Parens1{p|q,\mathcal O(\Hom0{\ger g/\ger h,\ger s/\ger t})}$ by $a(x_i)=\sum_{j=1}^n\xi_{ij}(a)y_j$ \fa $a\in\Hom0{\ger g/\ger h,\ger s/\ger t}$ (recall $V^*\subset\mathcal O(V)$ for $V$ any linear supermanifold). Then $\vphi^*(\xi)=(\vphi^*(\xi_{ij}))\in\GL\Parens1{p|q,\mathcal O(\mathcal U)}$ and 
	\[
	d=\Ber0{\vphi^*(\xi)}\in\mathcal O_{\mathcal G/\mathcal H}(U)\ .
	\]
\end{Lem}

\begin{Proof}
	In analogy with \cite[4.19]{schmitt-supergeom}, we may consider 
	\[
	\vphi\in\Parens1{\mathcal O_{\mathcal G/\mathcal H}(U)\otimes\Hom0{\ger g/\ger h,\ger s/\ger t}}_0\ .
	\]
	Then $\vphi(x_i)\in\mathcal O_{\mathcal G/\mathcal H}(U)\otimes\ger s/\ger t$ makes sense, is of parity $\Abs0{x_i}$, and we have the equation $\vphi(x_i)=\sum_{j=1}^n\vphi^*(\xi_{ij})y_j$. It follows that 
	\[
		D(\vphi(x_1),\dotsc,\vphi(x_n))=\Ber0{\vphi^*(\xi_{ij})}\cdot D(y_1,\dotsc,y_n)
	\]
	as elements of $\Ber[_{\mathcal O_{\mathcal G/\mathcal H}(U)}]0{\mathcal O_{\mathcal G/\mathcal H}(U)\otimes\ger s/\ger t}$. Similarly, the Berezinian sections may be considered as 
	\begin{gather*}
		\omega_{\mathcal G/\mathcal H}|_{\mathcal U}\in\Ber[_{\mathcal O_{\mathcal G/\mathcal H}(U)}]0{\mathcal O_{\mathcal G/\mathcal H}(U)\otimes(\ger g/\ger h)^*}\ ,\\
		\omega_{\mathcal S/\mathcal T}|_{\mathcal V}\in\Ber[_{\mathcal O_{\mathcal S/\mathcal T}(V)}]0{\mathcal O_{\mathcal S/\mathcal T}(V)\otimes(\ger s/\ger t)^*}\ .
	\end{gather*}

	For $\omega=\omega_{\mathcal G/\mathcal H}$, we have on $\mathcal U$:
	\[
	\tau^*\omega=\Ber0\tau^{-1}\circ\omega\circ\tau=\id\times\omega_o\equiv\omega_o\ ,
	\]
	by the $\mathcal G$-invariance of $\omega$ and the definition of $\Ber0\tau$ in \thmref{Par}{berchart}. A similar equation holds for $\omega_{\mathcal S/\mathcal T}$. Hence, abbreviating $D(x)=D(x_1,\dotsc,x_n)$, \emph{etc.}, we obtain
	\begin{align*}
		\Dual1{\tau^*\phi^*\omega_{\mathcal S/\mathcal T}}{D(x)}&=\Dual1{\sigma^*\omega_{\mathcal S/\mathcal T}}{D(\vphi(x))}\\
		&=\Ber0{\vphi^*(\xi_{ij})}\cdot\Dual1{(\omega_{\mathcal S/\mathcal T})_{o'}}{D(y)}\\
		&=\Ber0{\vphi^*(\xi_{ij})}\cdot\Dual1{(\omega_{\mathcal G/\mathcal H})_o}{D(x)}\\
		&=\Ber0{\vphi^*(\xi_{ij})}\cdot\Dual1{\tau^*\omega_{\mathcal G/\mathcal H}}{D(x)}
	\end{align*}
	where the first equality follows from \eqref{eq:locexpr}. This proves the assertion.
\end{Proof}

\begin{Prop}[prodsubsupergrp]
	Let $\mathcal U$ be a Lie supergroup and $\mathcal M$, $\mathcal H$ closed subsupergroups \scth the map $m=m_{\mathcal U}:\mathcal M\times\mathcal H\to\mathcal U$ is an isomorphism onto an open subsupermanifold $\mathcal V$ of $\mathcal U$. For a suitable normalisation of Berezinian sections, 
	\[
		\int_{\mathcal U}f=\int_{\mathcal M\times\mathcal H} m^*f\cdot\frac{\pr_2^*\Ber0{\Ad_{\ger h}}}{\pr_2^*\Ber0{\Ad_{\ger u}}}
	\]
	\fa $f\in\Gamma_c(U,\mathcal O_{\mathcal U})$ \scth $\supp f\subset V$. 
\end{Prop}

\begin{Proof}
	We consider $\mathcal U$ and $\mathcal M\times\mathcal H$ as quotients by the trivial subsupergroup. A standard calculation (for instance, using generalised points) shows that the local expression $\vphi:\mathcal M\times\mathcal H\to\Hom0{\ger m\oplus\ger h,\ger u}_0$ of $Tm$ is given by 
	\[
		\vphi=\iota_{\ger h\to\ger u}+\Ad^{\ger m\to\ger u}_{\mathcal H}\circ\, i_{\mathcal H}\circ\pr_2 
	\]
	where $\iota_{\ger h\to\ger u}$ is the constant morphism $\mathcal M\times\mathcal H\to\Hom0{\ger h,\ger u}$ whose value is the inclusion $\ger h\to\ger u$, and where $\Ad^{\ger m\to\ger u}_{\mathcal H}$ is the morphism $\mathcal H\to\Hom0{\ger m,\ger u}_0$ induced by restriction and corestriction by the adjoint morphism of $\mathcal U$. 
	
	It follows that $\vphi^*(\xi_{ij})=(\pr_2^*\Ber0{\Ad_{\ger h}})\cdot(\pr_2^*\Ber0{\Ad_{\ger u}})^{-1}$, so that the assertion follows from \thmref{Lem}{berhomsplociso} and the invariance of the Berezin integral under oriented isomorphisms. 
\end{Proof}

\end{document}